\documentclass[10pt,letterpaper]{article}
\usepackage[utf8]{inputenc}
\usepackage{authblk}
\usepackage{setspace}
\usepackage[margin=1.25in]{geometry}
\usepackage{graphicx}
\graphicspath{ {./figures/} }
\usepackage{subcaption}
\usepackage{commath}

\usepackage{bm}
\usepackage{enumitem}
\usepackage{algorithm}
\usepackage{amsmath,amssymb,mathtools,booktabs}

\usepackage{algpseudocode}

\usepackage{tabularx}
\usepackage{array}
\newcolumntype{Y}{>{\raggedright\arraybackslash}X}

\usepackage{float}
\usepackage{caption}

\usepackage{longtable}
\usepackage{xurl}
\newcolumntype{L}[1]{>{\raggedright\arraybackslash}p{#1}}
\newcommand{\feat}[1]{\begingroup\urlstyle{same}\path{#1}\endgroup}

\usepackage[backend=biber,
            style=authoryear-comp,
            sorting=nyt,
            natbib=true,
            maxcitenames=2,
            maxbibnames=99,
            uniquelist=false]{biblatex}
\usepackage{hyperref}
\hypersetup{hidelinks,
  pdftitle={Assortment Control Unlocks the Value of Dynamic Pricing in Mixed Last-Mile Delivery},
  pdfauthor={Ze Zhou and Balázs Kulcsár}}
\usepackage{cleveref}

\title{Assortment Control Unlocks the Value of Dynamic Pricing in Mixed Last-Mile Delivery}
\author[1]{Ze Zhou}
\author[1]{Balázs Kulcsár}

\affil[1]{Department of Electrical Engineering, Chalmers University of Technology, Gothenburg, Sweden.}

\date{}

\begin{document}

\maketitle

\begin{abstract}The growth of e-commerce has intensified the need for last-mile delivery systems that can jointly manage customer choice and operational efficiency. We study the Dynamic Offering and Pricing of Mixed Delivery Options (DOPMDO) problem, in which a logistics service provider dynamically selects and prices attended home-delivery time slots and out-of-home pickup options for sequentially arriving customers. Each decision affects immediate revenue, customer acceptance, and the route-dependent fulfillment cost realized at the end of the booking horizon. We formulate DOPMDO as a finite-horizon Markov decision process and propose State-Value Anchored Pricing via Approximate Dynamic Programming (SVAP--ADP). The method learns a continuation-value approximation on an aggregate state representation of the mixed-delivery system and uses accepted-versus-rejected value differences to estimate option-level opportunity costs. These opportunity costs are embedded in an anchored trust-region pricing problem around a calibrated fine-static benchmark. Computational experiments on the real-world Seattle instance show that SVAP--ADP increases mean episode profit by 7.0\% relative to current practice (95\% CI: 6.7--7.2\%), primarily by reducing terminal fulfillment cost while maintaining a stable home--locker--opt-out mix. Assortment-control experiments show that dynamic pricing is most effective when the menu exposes operationally valuable locker alternatives, with richer candidate pools delivering substantially larger gains than restricted nearest-locker menus. These results indicate that anticipatory pricing and assortment control are complementary: pricing steers customers toward lower-cost options, but the menu determines whether high-value consolidation opportunities are available in the first place. Ablation and sensitivity analyses further show that stable algorithm performance requires forward-looking opportunity-cost estimation, a sufficiently rich aggregate state representation, and anchored price control.
\end{abstract}

\section{Introduction}

Logistics service providers (LSPs) face a persistent profitability challenge in last-mile delivery. Although parcel volumes continue to grow, thin profit margins are often eroded by the high operating cost of attended home delivery (AHD), which requires vehicles to visit individual residences within promised service windows. Long handling times, urban congestion, parking constraints, and failed delivery attempts further increase the cost burden of AHD operations \citep{dalla2020commercial,ranjbari2023parcel,chen2017parking}. With global parcel volumes increasing from 64 billion in 2016 to more than 161 billion in 2022 and expected to reach 225 billion by 2028 \citep{pitneybowes2023}, these inefficiencies might translate into mounting financial and environmental pressure on logistics operations.

Out-of-home (OOH) delivery has gradually become an important operational alternative to mitigate these inefficiencies. Instead of serving each parcel at the customer's residence, LSPs consolidate deliveries at shared pickup locations such as parcel lockers, retail outlets, or grocery stores. The adoption of parcel lockers has expanded rapidly across Europe in recent years \citep{pinchasik2025replacing}; in Sweden, delivery points and parcel lockers account for 75\% of parcel deliveries \citep{vakulenko2023parcel}. By aggregating multiple parcels at common pickup points, OOH delivery can reduce driving distances \citep{enthoven2020two}, shorten service times, and lower the risk of failed delivery \citep{akkerman2025learning, savelsbergh201650th}, while remaining attractive to customers who value flexible pickup times.

Despite these advantages, operating AHD and OOH delivery jointly is far from straightforward. Customers differ in their preferences for delivery modes, time windows, prices, and pickup distances. Moreover, customer requests arrive sequentially during a booking horizon, while the final routing cost is realized only after all accepted orders have been confirmed. Each accepted order changes the spatial and temporal structure of the fulfillment problem and therefore changes the marginal value of future orders. A home-delivery acceptance may increase route fragmentation, while a locker acceptance may improve consolidation but require a discount to attract the customer. Conversely, overly aggressive surcharges or discounts may either induce excessive opt-outs or erode revenue. The LSP must therefore make real-time offering and pricing decisions that balance immediate revenue, customer choice, fleet capacity, and long-horizon route consolidation.

Dynamic offering and pricing have emerged as promising control levers for steering demand toward delivery options that are less costly to fulfil \citep{yang2017approximate,strauss2021dynamic,galiullina2024demand,akkerman2025learning}. By selecting a tailored menu of AHD time slots and OOH pickup options, and by attaching option-specific surcharges or discounts, an LSP can influence the distribution of accepted demand across delivery channels. However, designing effective policies for mixed AHD--OOH systems raises three methodological challenges. First, the system state contains all accepted orders, their selected delivery modes, time-window commitments, and tentative route structure, making exact dynamic programming computationally intractable. Second, the economic value of an accepted order depends not only on its immediate revenue but also on its effect on future capacity and route consolidation. Third, learned or myopic price adjustments can be unstable: if they depart too far from a well-calibrated static policy, they may sacrifice revenue, induce excessive opt-outs, or create fragmented route plans that are costly to serve. Existing studies have made important progress on dynamic pricing for home delivery, time-slot management, and OOH delivery, but they typically simplify at least one of the key dimensions: the joint assortment of AHD and OOH options, the sequential nature of customer arrivals, or the route-dependent fulfillment cost induced by mixed delivery choices.

In this paper, we study the Dynamic Offering and Pricing of Mixed Delivery Options (DOPMDO) problem. Customers arrive sequentially, and the LSP must decide which AHD time slots and OOH pickup options to offer and how to price them before the final fulfillment routes are known. Each decision therefore affects immediate revenue, customer choice, and the terminal routing cost of the accepted orders. We formulate this problem as a finite-horizon Markov Decision Process (MDP) and propose State-Value Anchored Pricing via Approximate Dynamic Programming (SVAP--ADP). The method estimates aggregate post-decision state value, translates accepted-versus-rejected value differences into option-level opportunity costs, and applies these costs within an anchored trust-region pricing problem. In this way, SVAP--ADP uses dynamic prices to exploit consolidation opportunities while keeping customer-facing prices close to a reliable static benchmark.

The main contributions are as follows.
\begin{enumerate}
    \item We introduce the DOPMDO problem as a sequential decision model for mixed attended home delivery and out-of-home pickup. The formulation integrates assortment selection, option-level pricing, stochastic customer choice, and route-dependent terminal fulfillment cost. It captures the key operational trade-off that an accepted order generates immediate revenue but also reshapes future consolidation and routing opportunities.
    \item We develop SVAP--ADP, an anchored value-based pricing method for large-scale mixed delivery systems. The method approximates aggregate post-decision state value, converts accepted-versus-rejected value differences into option-level opportunity costs, and combines these costs with immediate routing information inside a trust-region pricing problem. The resulting policy is forward-looking while remaining stable and interpretable.
    \item We evaluate SVAP--ADP on a real Seattle last-mile delivery network with both home-delivery and parcel-locker options. The experiments compare against a current-practice Static-uniform policy, a calibrated Static-fine benchmark, a myopic routing-aware policy, a compact value-function variant, and an unanchored FreePrice ablation. SVAP--ADP improves profit over both static baselines primarily by reducing terminal fulfillment cost, rather than by increasing revenue or suppressing demand.
    \item We identify the mechanisms that make dynamic pricing effective in mixed AHD--OOH systems. The results show that immediate insertion costs are too local to capture future consolidation value, unrestricted learned prices can destabilize the demand mix, and a richer aggregate state representation is needed for reliable opportunity-cost pricing. The sensitivity analyses further show that pricing and assortment control are complementary: prices steer demand, while the candidate pool determines whether operationally useful locker options are available.
\end{enumerate}

The remainder of this paper is organized as follows. Section~\ref{sec:literature} reviews related work on dynamic offering, pricing, AHD and OOH delivery. Section~\ref{sec:problem} introduces the DOPMDO problem and its MDP formulation. Section~\ref{sec:method} presents the SVAP--ADP method. Section~\ref{sec:seattle} reports the numerical study on the Seattle network. Section~\ref{sec:conclusion} concludes the paper and discusses directions for future research.

\section{Literature review}
\label{sec:literature}

In this section, we first review demand management for AHD, followed by studies on OOH delivery and mixed-channel delivery design. We then discuss sequential decision methods for integrated demand management and vehicle routing problems, and position our paper at the intersection of these research streams.

\subsection{Demand management in attended home delivery}
Demand management for AHD studies how an LSP can influence customer choices before the final delivery routes are constructed. Early work showed that delivery-slot availability, incentives, and acceptance decisions can be used to steer demand toward operationally attractive regions and time windows. \citet{campbell2006incentive} study dynamic acceptance and incentive schemes for home delivery services, highlighting the operational value of influencing customer time-window choices. \citet{agatz2011time} consider time-slot management at the area level and use continuous approximations to anticipate delivery costs. Studies by \citet{ehmke2014customer} and \citet{visser2019strategic} further demonstrate that controlling slot availability can improve routing feasibility and reduce fulfillment cost. Recent reviews by \citet{wassmuth2023demand} and \citet{fleckenstein2023recent} provide broader classifications of AHD demand management and integrated demand management and vehicle routing problems.

A second branch of this literature studies pricing rather than pure availability control. \citet{asdemir2009dynamic} formulate dynamic time-slot pricing with a multinomial logit choice model, but rely on simplified capacity representations. \citet{yang2016choice} combine a choice-based pricing model with dynamically estimated delivery costs, using real e-grocery data to show that customers can be steered through delivery charges. \citet{yang2017approximate} extend this line by using approximate dynamic programming (ADP) to account for future revenue and routing-cost effects. \citet{klein2019differentiated} study differentiated static time-slot pricing under routing considerations, while \citet{koch2020route} develop a route-based ADP approach for dynamic pricing in AHD. \citet{vinsensius2020dynamic} similarly study dynamic incentives for slot management in e-commerce AHD. More recently, \citet{abdollahi2023demand} incorporate forecast orders into dynamic routing to support time-slot demand management. \citet{abdolhamidi2025tactical} further incorporate heterogeneous customer preferences into the design of slot assortments and price discounts. Their numerical results show that using a mixed logit model enables the LSP to tailor its service offerings more effectively, thereby increasing profitability and improving customer retention.

Several papers relax the classical assumption that customers choose exactly one narrow time window. \citet{yildiz2020pricing} study discounts for delivery-time flexibility across multiple periods. \citet{strauss2021dynamic} introduce flexible time slots, where customers accept uncertainty about the final delivery window in exchange for a lower delivery charge, and develop a dynamic pricing policy informed by approximate opportunity costs. These studies demonstrate the value of using pricing to obtain operational flexibility. However, they focus mainly on temporal flexibility in home delivery. They do not jointly decide prices and assortments over home-delivery time slots and spatially distinct out-of-home pickup locations.

\subsection{Out-of-home delivery and mixed delivery-option design}
Out-of-home (OOH) delivery has received increasing attention as a way to consolidate parcels, reduce failed deliveries, and improve last-mile efficiency. Empirical and optimization studies show that parcel lockers and pickup points can reduce delivery effort and create operational economies of density. \citet{ranjbari2023parcel} provide field evidence on the effect of parcel lockers on delivery times. \citet{enthoven2020two} study a two-echelon vehicle routing problem with covering options, in which parcel lockers serve as shared delivery locations. \citet{mancini2021vehicle} introduce vehicle routing with private and shared delivery locations, combining home delivery with alternative delivery points. \citet{lin2022profit} and \citet{lyu2022last} study parcel-locker location and locker alliance network design, while \citet{mancini2023locker} address locker-location planning under uncertainty in demand and capacity. \citet{janinhoff2024out} provide a recent review of OOH delivery optimization, covering facility location, routing, location-routing, and emerging operational challenges.

Most OOH studies, however, are static: customer locations, potential delivery options, or aggregate demand are known before decisions are made. Some works incorporate customer preferences or incentives, but typically without modeling a fully sequential booking process. \citet{dumez2021large} investigate a setting in which customers can actively specify delivery options, each associated with time windows and preference levels. They formulate the problem as a mixed-integer linear program and develop a large neighborhood search method to solve it. Their results show that the proposed formulation can generate substantial cost savings while maintaining a high quality of service. \citet{zhang2023joint} study joint location and pricing optimization for self-service delivery under customer choice, with a focus on strategic facility and price design. \citet{galiullina2024demand} examine demand steering with home and pickup-point delivery options, integrating incentive and routing decisions under uncertain customer acceptance. \citet{zhou2026last} investigate a mixed last-mile delivery problem with heterogeneous customer behavior and formulate it as a two-stage stochastic program, showing that joint assortment and pricing decisions can substantially improve delivery performance. Together, these studies demonstrate the operational value of steering customers toward pickup points and jointly managing home and OOH delivery options. However, they all assume that the relevant customer demand is known before decisions are made, and therefore do not address the real-time problem of jointly selecting and pricing both home time slots and nearby OOH locations as customers arrive sequentially.

The most closely related OOH study is \citet{akkerman2025learning}, who introduce dynamic selection and pricing of OOH delivery options. They formulate a sequential decision problem in which each arriving customer can be offered a subset of OOH locations with discounts or charges, and they propose a machine-learning policy using a spatial-temporal state representation. Their work is important because it moves OOH delivery from static optimization to sequential offering and pricing. Nevertheless, their setting focuses on choosing between home delivery and OOH locations without modeling attended home-delivery time slots. In contrast, our problem considers a mixed menu consisting of multiple AHD time slots and multiple nearby lockers, so the provider must jointly manage both the temporal dimension of AHD and the spatial consolidation dimension of OOH delivery.

\subsection{Sequential decision for last-mile delivery}
Because customers arrive dynamically and each accepted order reshapes both the routing problem and future acceptance opportunities, many last-mile pricing problems are naturally modeled as finite-horizon Markov decision processes. Exact dynamic programming is generally intractable at realistic scales, so the literature relies on decomposition, opportunity-cost approximation, ADP, simulation-based policies, and, more recently, machine learning. \citet{fleckenstein2023recent} review integrated demand management and vehicle routing problems and emphasize the common structure in which providers control prices, availability, or acceptance decisions during a booking horizon while fulfillment is performed by a vehicle fleet. \citet{fleckenstein2025concept} further analyze the concept of opportunity cost in integrated demand management and vehicle routing, showing why opportunity-cost approximation is more difficult than in classical revenue management: the value of a request depends jointly on routing interactions, capacity consumption, and displaced future revenue.

ADP-based pricing has been particularly influential in AHD and same-day delivery. \citet{yang2017approximate} use ADP to estimate the future value of capacity and routing resources for time-slot pricing. \citet{koch2020route} construct route-based value approximations for dynamic AHD pricing. \citet{ulmer2020dynamic} studies dynamic pricing and routing for same-day delivery, where booking and service horizons overlap, and proposes an anticipatory pricing and routing policy. \citet{klein2023dynamic} combine dynamic demand management with online tour planning for same-day delivery, optimizing delivery spans and prices for incoming requests. \citet{banerjee2025pricing} consider pricing and demand management for integrated same-day and next-day delivery systems, further illustrating the growing interest in coordinating pricing with fleet operations. These studies confirm the value of anticipatory pricing, but most of them control temporal delivery promises rather than mixed AHD--OOH delivery options. Machine-learning methods have also begun to appear in last-mile demand management. \citet{akkerman2025learning} use a convolutional neural network with a spatial-temporal state encoding for dynamic OOH selection and pricing. These methods are attractive because they can capture complex spatial and temporal interactions. At the same time, they may be unstable in revenue-management settings if learned adjustments deviate too far from a reliable benchmark policy, especially when value estimates are noisy.

\begin{table}[t]
\centering
\caption{Comparison of related research on last-mile demand management.}
\label{tab:literature_classification}
\footnotesize
\setlength{\tabcolsep}{4pt}
\renewcommand{\arraystretch}{1.08}

\begin{tabular}{@{}lcccc@{}}
\toprule
\textbf{Article} & \textbf{Setting} & \textbf{Decision} & 
\textbf{Approach and choice model} & \textbf{Channel} \\
\midrule

\multicolumn{5}{@{}l}{\emph{AHD: time-slot offering and pricing}} \\
\citet{yang2016choice} & D & P & Insertion-cost approx. (MNL) & AHD slots \\
\citet{yang2017approximate} & D & P & ADP (MNL) & AHD slots \\
\citet{klein2019differentiated} & S & P & MIP+ADP (Nonparametric) & AHD slots \\
\citet{koch2020route} & D & P & Route-based ADP (MNL) & AHD slots \\
\citet{strauss2021dynamic} & D & P & LP-based opportunity cost approx. (MNL) & AHD slots \\
\citet{abdollahi2023demand} & D & P & Forecast rollout (MNL) & AHD slots \\
\citet{abdolhamidi2025tactical} & S & O+P & MIP (Mixed Logit) & AHD slots \\

\addlinespace[1pt]
\midrule
\multicolumn{5}{@{}l}{\emph{Same-day delivery}} \\
\citet{ulmer2020dynamic} & D & P & Meso-level ADP (MNL) & AHD slots \\
\citet{klein2023dynamic} & D & O+P & ADP + online tour planning (MNL) & AHD slots \\

\addlinespace[1pt]
\midrule
\multicolumn{5}{@{}l}{\emph{OOH and mixed-channel delivery}} \\
\citet{enthoven2020two} & S & O & Two-echelon routing optimization (--) & OOH \\
\citet{dumez2021large} & S & O & MIP (Preference-level) & AHD+OOH \\
\citet{mancini2021vehicle} & S & O & MIP (--) & AHD+OOH \\
\citet{galiullina2024demand} & S & P & MIP (MNL) & OOH \\
\citet{zhang2023joint} & S & O+P & Choice-based optimization (MNL) & OOH \\
\citet{akkerman2025learning} & D & O+P & CNN-based learning policy (MNL) & OOH \\
\citet{zhou2026last} & S & O+P & MIP (MNL) & AHD+OOH \\

\midrule
\textbf{This paper} & \textbf{D} & \textbf{O+P} & 
\textbf{SVAP--ADP (MNL)} & \textbf{AHD+OOH} \\

\bottomrule
\end{tabular}

\vspace{1mm}
\begin{minipage}{0.99\textwidth}
\scriptsize
\emph{Note.}
D = dynamic; S = static; O = option, assortment, location, or assignment decision; P = pricing or incentive decision; O+P = joint option and pricing decision.
AHD = attended home delivery; OOH = out-of-home delivery; ADP = approximate dynamic programming; CNN = convolutional neural network; LP = linear programming; MIP = mixed-integer programming; MNL = multinomial logit.
\end{minipage}
\end{table}

To position our study relative to the most closely related literature, Table~\ref{tab:literature_classification} classifies representative papers by decision setting, demand-management lever, methodology, choice model, and delivery channel. The table highlights that existing dynamic offering and pricing studies typically focus either on AHD time slots or on OOH delivery locations, while studies that combine home and OOH delivery are mostly static or do not jointly control both assortment and prices in an online setting.

Our paper addresses this gap by introducing the DOPMDO problem as a sequential decision model for mixed last-mile delivery. Unlike prior studies that focus primarily on AHD time-slot control, OOH selection and pricing, or static mixed-channel design, we jointly consider dynamic pricing and assortment control for both attended home-delivery time slots and out-of-home pickup locations as customers arrive, while accounting for stochastic choice behavior and route-dependent terminal fulfillment cost.

\section{Problem description and model}
\label{sec:problem}

We study the DOPMDO problem faced by a last-mile delivery provider during a booking horizon before route execution. Customers arrive sequentially. At each arrival, the provider observes the customer's location and feasible delivery alternatives, then chooses which options to display and what price adjustment to attach to each displayed option. The option set combines AHD in a finite set of delivery time windows and OOH delivery to parcel lockers. The customer's response is stochastic. Accepted orders are fulfilled only after the booking horizon ends, so each online decision affects both immediate revenue and the route-dependent fulfillment cost realized at the terminal stage.

\subsection{Problem narrative}
\label{subsec:problem-narrative}

The provider operates a fleet $\mathcal V=\{1,\ldots,V\}$ of homogeneous vehicles with capacity $Q$ from a single depot. The service region is represented by a travel-time network $G=(\mathcal N^0,\mathcal E)$, where $\mathcal N^0$ contains the depot, customer home locations revealed during the booking horizon, and a fixed set $\mathcal L$ of parcel lockers. Travel time between locations $i$ and $j$ is denoted by $\tau(i,j)$. The set of AHD time windows is denoted by $\mathcal S$. A home-delivery option assigns the current order to the customer's home location in one of these time windows. A locker option assigns the parcel to one of the candidate lockers considered for the current customer.

The booking horizon contains a finite sequence of customer-arrival opportunities. We index decision epochs by $k$ and let $D$ denote the maximum horizon length used in the dynamic program and in the value approximation. The realized arrival process may terminate before $D$; this is represented by a terminal realization of the exogenous information. At epoch $k$, the current customer is denoted by $C_k^{\mathrm{new}}$. The observed customer information consists of the home location $h_k$, the service or load requirement $b_k$, and a finite candidate locker set $\mathcal L_k\subseteq\mathcal L$. The candidate locker set may be generated from nearby lockers, active lockers, or another customer-relevance rule specified by the online policy; the MDP only requires it to be finite and known before the decision.

The provider controls two customer-facing levers. First, it chooses the displayed assortment, namely the subset of feasible AHD and OOH options shown to the customer. Second, it chooses an option-specific price adjustment. Price adjustments are measured relative to the base order revenue $r$: positive values are surcharges and negative values are discounts. The customer chooses from the displayed menu according to a stochastic discrete-choice model with an outside option.

The economic trade-off is dynamic. AHD orders generate revenue but require residential visits within selected time windows, which may fragment routes across both space and time. OOH orders can consolidate multiple parcels at shared pickup locations, but discounts used to stimulate OOH demand reduce immediate revenue. Opening many weakly used lockers may also increase route stops, while reusing an already active locker may add little or no additional travel. Because the final routes are constructed only after all accepted orders have been confirmed, the provider must price and offer current options based on both current revenue and the future operational flexibility left for subsequent customers.

\subsection{Markov decision process}
\label{subsec:mdp}

We formulate DOPMDO as a finite-horizon MDP. The decision epoch begins after the current customer has arrived and revealed the information needed for offering and pricing. Let
\begin{equation}
B_k=(k,\Omega_k,P_k)
\label{eq:booking-state}
\end{equation}
denote the booking state inherited from previous decisions, where $k$ is the epoch index, $\Omega_k$ is the set of accepted orders, and $P_k$ is the tentative fulfillment plan used for feasibility checks and successor-state construction. The full decision state is
\begin{equation}
S_k=(B_k,C_k^{\mathrm{new}}),
\label{eq:decision-state}
\end{equation}
where $C_k^{\mathrm{new}}=(h_k,b_k,\mathcal L_k)$ is the current customer information. Thus, $S_k$ is the state on which the online decision is made, while $B_k$ is the post-response booking state carried forward from the previous epoch.

The feasible option set is generated from $S_k$. AHD options are $(h_k,w)$ for $w\in\mathcal S$, and OOH options are lockers $\ell\in\mathcal L_k$. An option is feasible if the insertion operator can construct a feasible option-contingent plan. For home delivery, the plan must respect the selected time window and vehicle-capacity restrictions; for locker delivery, it must respect vehicle capacity and, when modeled, residual locker capacity. If the route plan already visits locker $\ell$, an additional parcel assigned to $\ell$ can be consolidated at that stop.

A decision at epoch $k$ is denoted by $x_k=(U_k,p_k,\{P_k^o\}_{o\in U_k})$. Here, $U_k$ is the displayed option set, $p_k(o)$ is the price adjustment for option $o$, and $P_k^o$ is the option-contingent plan obtained if option $o$ is accepted. The feasible action set is denoted by $\mathcal A(S_k)$. The empty menu is included in $\mathcal A(S_k)$ and represents a provider-controlled rejection action: if $U_k=\emptyset$, then $q_k(0\mid S_k,x_k)=1$ and the system follows the opt-out successor $B^0_{k+1}$.

The customer's choice follows a multinomial logit model. The systematic utility of option $o$ is
\begin{equation}
v_k(o)=
\begin{cases}
\mu_w+\beta^p p_k(o), & o=(h_k,w),\\
\mu_\ell-\beta^d d(h_k,\ell)+\beta^p p_k(o), & o=\ell\in\mathcal L_k,
\end{cases}
\label{eq:utility}
\end{equation}
where $\mu_w$ is the baseline utility of AHD slot $w$, $\mu_\ell$ is the OOH baseline utility, $\beta^p<0$ is price sensitivity, $\beta^d>0$ is locker-distance sensitivity, and $d(h_k,\ell)$ is the distance from the customer's home to locker $\ell$. Normalizing the outside-option utility to zero, the choice probabilities are

\begin{equation}
q_k(u\mid S_k,x_k)
=
\frac{\exp(v_k(u))}
{\sum_{\bar u\in U_k\cup\{0\}}\exp(v_k(\bar u))},
\qquad u\in U_k\cup\{0\}.
\label{eq:mnl_choice_probability}
\end{equation}
Although the choice probabilities depend on $x_k$ only through the displayed
menu $U_k$ and prices $p_k$, conditioning on $(S_k,x_k)$ keeps the notation
aligned with the Bellman recursion and the option-contingent successor states.
If $U_k=\emptyset$, then $q_k(0\mid S_k,x_k)=1$.

Let $u_k\in U_k\cup\{0\}$ denote the realized choice, where $0$ is the outside option. If $u_k=0$, no order is appended. If $u_k=o\in U_k$, the current order is added to the accepted-order set and the tentative plan becomes $P_k^o$. The corresponding booking-state successors are
\begin{equation}
B_{k+1}^{0}=(k+1,\Omega_k,P_k),
\qquad
B_{k+1}^{o}=\bigl(k+1,\Omega_k\cup\{(C_k^{\mathrm{new}},o)\},P_k^o\bigr).
\label{eq:successor-booking-states}
\end{equation}
The next decision state is formed by combining the realized successor booking state with the next customer request, unless the booking horizon has ended.

At the terminal stage, the provider executes routes to serve all accepted orders. The realized fulfillment cost is recomputed from the terminal booking state and decomposes into three cost categories:
\begin{equation}
C^{\mathrm{ful}}(B_\tau)
=
\gamma\left(
C_\tau^{\mathrm{dist}}
+
C_\tau^{\mathrm{time}}
+
C_\tau^{\mathrm{fix}}
\right),
\label{eq:fulfillment_cost}
\end{equation}
where $C_\tau^{\mathrm{dist}}$ is distance-dependent routing cost, $C_\tau^{\mathrm{time}}$ is service time at home and locker stops, $C_\tau^{\mathrm{fix}}$ is fixed route and vehicle cost, and $\gamma$ is a fulfillment-cost multiplier. 

The provider seeks a policy $\pi$ that maximizes expected profit,
\begin{equation}
\max_{\pi\in\Pi}
\mathbb E^\pi
\left[
\sum_{k=1}^{\tau-1}
\bigl(r+p_k(u_k)\bigr)\mathbf 1\{u_k\neq 0\}
-
\pi_{out} \mathbf 1\{u_k= 0\}
-
C^{\mathrm{ful}}(B_\tau)
\right],
\label{eq:objective}
\end{equation}
where $p_k(0)=0$ by convention and $\tau$ denotes the terminal epoch. Thus, revenue $r+p_k(u_k)$ is earned only when the customer accepts an offered option, an opt-out incurs penalty $\pi_{out}$, and fulfillment cost is incurred once at the terminal booking state.

\subsection{Decision-state and post-decision value functions}
\label{subsec:bellman}

Let $V_k(S_k)$ denote the optimal expected profit from epoch $k$ onward after the current customer has been observed. The value is defined on the full decision state $S_k=(B_k,C_k^{\mathrm{new}})$ because the current customer determines the feasible options, choice probabilities, and price adjustments. The exact Bellman recursion is
\begin{equation}
V_k(S_k)
=
\max_{x_k\in\mathcal A(S_k)}
\sum_{u\in U_k\cup\{0\}}
q_k(u\mid S_k,x_k)
\left[
R_k(u;x_k)
+
V_{k+1}^{x}\!\left(B_{k+1}^{u}\right)
\right],
\label{eq:bellman-full-state}
\end{equation}
where $q_k(u\mid S_k,x_k)$ is the MNL choice probability, $R_k(0;x_k)=-\pi_{out}$ for opt-out, and $R_k(o;x_k)=r+p_k(o)$ for $o\in U_k$.

The post-decision value $V^x_{k+1}(B)$ is the conditional expected continuation value after the current customer's outcome has been incorporated, given the resulting booking state $B$, but before the next customer is observed:
\begin{equation}
V^x_{k+1}(B)
=
\mathbb{E}\!\left[
V_{k+1}(B,\widetilde C_{k+1})
\,\middle|\, B
\right],
\label{eq:post-decision-value}
\end{equation}

where $\widetilde C_{k+1}$ is either the next customer request or the terminal signal $\emptyset$. We use the terminal convention
\begin{equation}
V_{k+1}(B,\emptyset)=-C^{\mathrm{ful}}(B),
\label{eq:terminal_convention}
\end{equation}
so the route-dependent terminal cost enters the recursion through the terminal realization.

The timing within one epoch can be summarized as
\begin{equation}
B_k
\xrightarrow{\;C_k^{\mathrm{new}}\;}
S_k=(B_k,C_k^{\mathrm{new}})
\xrightarrow{\;x_k=(U_k,p_k)\;}
u_k
\xrightarrow{\;\text{update }\Omega_k\text{ and }P_k\;}
B_{k+1}^{u_k}
\xrightarrow{\;\widetilde C_{k+1}\;}
S_{k+1}=(B_{k+1}^{u_k},\widetilde C_{k+1}).
\label{eq:state-sequence}
\end{equation}
This sequence separates the transient customer request from the persistent booking state. The current customer enters the decision state only to define the feasible options, choice probabilities, and option-contingent route updates. After the response, the system carries forward only the updated booking state $B_{k+1}^{u_k}$, and the next decision state is formed when the next request is realized.

The MDP formulation makes clear why exact dynamic programming is not computationally viable: the booking state contains both the accepted-order set and the tentative fulfillment plan, so the state space grows rapidly over the booking horizon. This motivates the solution approach developed next. Rather than evaluating the full Bellman recursion, SVAP--ADP approximates the continuation value of post-response booking states using aggregate summaries of $B_k=(k,\Omega_k,P_k)$. Each feasible option is then assessed by comparing the successor state created by accepting that option with the opt-out successor, producing an option-level opportunity-cost signal for the online pricing problem.

\section{Solution approach: State-Value Anchored Pricing via ADP}
\label{sec:method}

The Bellman equation \eqref{eq:bellman-full-state} shows why online offering and pricing require anticipation. If the continuation value of each successor booking state were known, the provider could evaluate the future consequence of accepting each feasible option and solve a one-step menu-pricing problem at every arrival. This is the logic used in dynamic time-slot pricing and same-day pricing-and-routing models: the dynamic program is too large to solve exactly, but value differences quantify the future cost of accepting a current order. In DOPMDO, these value differences must account for both AHD time-window commitments and OOH locker consolidation.

We propose State-Value Anchored Pricing via Approximate Dynamic Programming (SVAP--ADP). The method has three design elements. First, it learns an aggregate approximation of the post-decision continuation value from simulated booking trajectories. Second, it estimates option-level opportunity costs by comparing the value of the opt-out successor state with the value of the accepted successor state. Third, it converts these opportunity costs into customer-facing prices within a trust region around a calibrated fine-static price anchor. The value approximation provides anticipation; the anchor and trust region stabilize the price optimization and prevent extreme reactions to approximation error.

\subsection{Opportunity costs and the one-step pricing problem}
\label{subsec:opportunity-cost}

The post-decision value in \eqref{eq:post-decision-value} allows the Bellman recursion to be written in the same opportunity-cost form used in dynamic pricing for delivery time slots. An accepted option earns immediate revenue, whereas an opt-out incurs the immediate penalty $\pi_{\mathrm{out}}\ge 0$. Substituting \eqref{eq:post-decision-value} into \eqref{eq:bellman-full-state} gives
\begin{equation}
\begin{aligned}
V_k(S_k)=\max_{x_k\in\mathcal A(S_k)}\Bigg\{&
q_k(0\mid S_k,x_k)\left[-\pi_{\mathrm{out}}+V_{k+1}^{x}(B_{k+1}^{0})\right] \\
&+\sum_{o\in U_k}q_k(o\mid S_k,x_k)\left[r+p_k(o)+V_{k+1}^{x}(B_{k+1}^{o})\right]\Bigg\}.
\end{aligned}
\label{eq:bellman-post-decision}
\end{equation}

For an offered option $o\in U_k$, define the exact \textit{opportunity cost} of acceptance as the loss in conditional post-decision continuation value relative to the opt-out successor:
\begin{equation}
\Delta_k(S_k,o)=V_{k+1}^{x}(B_{k+1}^{0})-V_{k+1}^{x}(B_{k+1}^{o}).
\label{eq:exact-opp-cost}
\end{equation}
Because $V_{k+1}^{x}(\cdot)$ is a future-only post-decision value, $\Delta_k(S_k,o)$ measures the expected future value displaced by accepting the current customer under option $o$ rather than carrying forward the opt-out successor state. This displacement includes the expected change in terminal fulfillment cost, future revenue opportunities, future opt-out penalties, and operational flexibility. It does not include the current immediate revenue or the current opt-out penalty; these current-epoch terms enter the pricing objective explicitly below. Note that the superscript $x$ in $V_{k+1}^{x}$ denotes the post-decision value, not dependence on the current action.

The rearrangement follows by adding and subtracting the opt-out continuation value inside the accepted-option terms. For any fixed feasible action $x_k$, let $q_k^0=q_k(0\mid S_k,x_k)$, $q_k^o=q_k(o\mid S_k,x_k)$, $V^0=V_{k+1}^{x}(B_{k+1}^{0})$, and $V^o=V_{k+1}^{x}(B_{k+1}^{o})$. Since $q_k^0+\sum_{o\in U_k}q_k^o=1$,
\begin{equation}
\begin{aligned}
&q_k^0\left[-\pi_{\mathrm{out}}+V^0\right]+\sum_{o\in U_k}q_k^o\left[r+p_k(o)+V^o\right] \\
&\quad=V^0+\sum_{o\in U_k}q_k^o\left[r+p_k(o)+V^o-V^0\right]-q_k^0\pi_{\mathrm{out}} \\
&\quad=V^0+\sum_{o\in U_k}q_k^o\left[r+p_k(o)-\Delta_k(S_k,o)\right]-q_k^0\pi_{\mathrm{out}}.
\end{aligned}
\label{eq:opp-cost-rearrangement}
\end{equation}

Returning to the original notation, \eqref{eq:bellman-post-decision} becomes
\begin{equation}
V_k(S_k)=V_{k+1}^{x}(B_{k+1}^{0})+\max_{x_k\in\mathcal A(S_k)}
\left\{
\sum_{o\in U_k}q_k(o\mid S_k,x_k)\left[r+p_k(o)-\Delta_k(S_k,o)\right]
-q_k(0\mid S_k,x_k)\pi_{\mathrm{out}}
\right\}.
\label{eq:bellman-opportunity-cost}
\end{equation}
The first term is the conditional continuation value of the opt-out successor booking state and is common to all feasible actions. The maximized term contains the incremental expected contribution of the offered options and the expected opt-out penalty. Each accepted option contributes immediate revenue plus the price adjustment, net of the option's opportunity cost. The outside option contributes no accepted-order revenue and incurs penalty $\pi_{\mathrm{out}}$. The opt-out term must remain inside the maximization because $q_k(0\mid S_k,x_k)$ depends on the displayed menu and prices. Thus, if the exact opportunity costs were known, the online decision would solve the one-step pricing problem
\begin{equation}
x_k^*
\in
\arg\max_{x_k=(U_k,p_k)\in\mathcal A(S_k)}
\left\{
\sum_{o\in U_k}q_k(o\mid S_k,x_k)\left[r+p_k(o)-\Delta_k(S_k,o)\right]
-q_k(0\mid S_k,x_k)\pi_{\mathrm{out}}
\right\}.
\label{eq:pricing-exact-opp}
\end{equation}
The maximization is over the insertion-feasible action set $\mathcal A(S_k)$ defined in Section~\ref{subsec:mdp}: infeasible home time-window insertions, vehicle-capacity violations, and modeled locker-capacity violations are excluded before pricing, and the empty menu remains feasible as a provider-controlled rejection action. When $\pi_{\mathrm{out}}=0$, \eqref{eq:pricing-exact-opp} reduces to the no-penalty opportunity-cost pricing problem. Each accepted option earns immediate revenue $r+p_k(o)$ and consumes opportunity cost $\Delta_k(S_k,o)$. A low opportunity cost indicates that the option fits well with the current AHD route structure or locker consolidation pattern. A high opportunity cost indicates that acceptance is expected to reduce future revenue opportunities, consume scarce operational flexibility, or increase the terminal fulfillment cost. The opt-out penalty discourages menus and prices that reduce operating cost only by pushing too many customers to the outside option. At the final booking decision, the terminal convention in \eqref{eq:terminal_convention} implies
\[
\Delta_D(S_D,o)=C^{\mathrm{ful}}(B_{D+1}^{o})-C^{\mathrm{ful}}(B_{D+1}^{0}),
\]
so the last decision directly trades off immediate revenue, the expected opt-out penalty, and the change in terminal fulfillment cost. At earlier epochs, the same terminal cost and downstream demand effects affect pricing through the continuation values.

The exact opportunity cost in \eqref{eq:exact-opp-cost} is unavailable because $V^x$ is intractable. SVAP--ADP therefore learns an approximation $\widehat V^x$ on aggregate summaries of the post-decision booking state and computes
\begin{equation}
\widehat\Delta_k(S_k,o)=\widehat V^{x}(B_{k+1}^{0})-\widehat V^{x}(B_{k+1}^{o}).
\label{eq:approx-opportunity-cost}
\end{equation}
The learned object is the post-decision value function; the pricing input is the accepted-versus-opt-out value difference. Current-epoch revenue and the current opt-out penalty are not absorbed into $\widehat\Delta_k(S_k,o)$; they enter the one-step pricing objective explicitly. The remaining subsections describe the aggregate state representation, the offline value-learning procedure, and the anchored online pricing problem used to implement \eqref{eq:pricing-exact-opp} with approximate opportunity costs.

\subsection{Core aggregate variables and value representation}
\label{subsec:aggregate-state}

The value approximation targets the post-decision value $V^x$, not the full decision-state value $V_k(S_k)$. The aggregate representation is therefore derived from the booking state $B_k=(k,\Omega_k,P_k)$, after the previous customer outcome has been incorporated and before the next customer is observed. The current customer information $C_k^{\mathrm{new}}$ is excluded from the learned state representation because it affects the current feasible options, prices, and choice probabilities, but it does not persist after the customer either accepts an option or opts out.

The exact booking state is too detailed for direct value approximation because it contains all accepted commitments and the tentative fulfillment plan. We therefore use a compact set of operational summaries and let the implementation-level basis be generated as deterministic transformations of these summaries. Let $n_k$ be the number of accepted orders in $\Omega_k$, and let $n_k^H$ be the number of accepted AHD orders. The number of accepted OOH orders is derived as $n_k^L=n_k-n_k^H$. Thus, $(k,n_k,n_k^H)$ records booking progress, total committed load, and the realized AHD--OOH modal split.

To describe where the accepted commitments are concentrated, we use two active-resource maps. Partition the service area into coarse home-delivery cells and combine each cell with each AHD time window. Let $a=(g,w)$ denote a home cell--slot pair, and let $m_{ak}^H$ be the number of accepted AHD orders assigned to pair $a$ before epoch $k$. The active AHD map is
\begin{equation}
 A_k^H=\{a=(g,w):m_{ak}^H>0\},
\label{eq:active-home-map}
\end{equation}
understood together with the active loads $\{m_{ak}^H:a\in A_k^H\}$. Similarly, let $m_{\ell k}^L$ be the number of accepted OOH parcels assigned to locker $\ell$. The active locker map is
\begin{equation}
 A_k^L=\{\ell\in\mathcal L:m_{\ell k}^L>0\},
\label{eq:active-locker-map}
\end{equation}
understood together with the active locker loads $\{m_{\ell k}^L:\ell\in A_k^L\}$. Finally, let $C_k$ denote a raw route-burden summary computed from the tentative plan $P_k$, such as the accumulated travel distance used by the insertion routine. This quantity is a state summary used for learning; it is not the terminal fulfillment cost $C^{\mathrm{ful}}$ in the exact objective. Note that $P_k$ is constructed first, and $C_k$ is then evaluated as an aggregated feature or a metric of the current route plan. Table~\ref{tab:aggregate-dimensions} summarizes the aggregate variables and operational dimensions used to construct the SVAP--ADP value approximation.

The core aggregate information used by SVAP--ADP is
\begin{equation}
\mathcal X(B_k)=
\left(k,\; n_k,\; n_k^H,\; A_k^H,\; A_k^L,\; C_k\right).
\label{eq:core-aggregate-state}
\end{equation}
AHD route dispersion is derived from $A_k^H$ and the loads $m_{ak}^H$; locker consolidation and fragmentation are derived from $A_k^L$ and the loads $m_{\ell k}^L$; and modal balance is derived from $n_k$ and $n_k^H$. The value approximation is defined on a basis generated from $\mathcal X(B_k)$. Let $\phi(B_k)=\left(1,\phi_1(B_k),\ldots,\phi_p(B_k)\right)$ be the aggregate basis, where each non-intercept component is a deterministic transformation or normalization of the core variables in \eqref{eq:core-aggregate-state}. The full basis implemented is reported in Appendix~\ref{app:vfa_basis_terms}. 

\begin{table}[t]
\centering
\caption{Core aggregate variables and derived state dimensions used by SVAP--ADP.}
\label{tab:aggregate-dimensions}
\begin{tabularx}{\textwidth}{@{}L{0.22\textwidth}L{0.27\textwidth}Y@{}}
\toprule
Derived dimension & Core variable(s) & Interpretation in the value approximation \\
\midrule
Booking progress & $k$ & Position in the booking horizon; determines how much future demand and routing opportunity remain. \\
Committed load & $n_k$ & Total demand already committed to the fulfillment system. \\
Modal balance & $n_k^H$ ($n_k^L=n_k-n_k^H$) & Split of accepted demand between AHD and OOH delivery. \\
AHD route dispersion & $A_k^H$ and loads $m_{ak}^H$ & Concentration or fragmentation of home commitments across cell--slot route patterns. \\
OOH consolidation and fragmentation & $A_k^L$ and loads $m_{\ell k}^L$ & Reuse of active lockers, opening of new lockers, and singleton-locker fragmentation. \\
fulfillment burden & $C_k$  & Route-burden information not captured by counts and active-resource maps alone. \\
\bottomrule
\end{tabularx}
\end{table}

Let $z_j(B_k)$ denote the standardized value of basis component $\phi_j(B_k)$. The approximate post-decision value is represented by a linear model,
\begin{equation}
\widehat V^x(B_k)=\theta_0+\sum_{j=1}^{p}\theta_j z_j(B_k).
\label{eq:vfa}
\end{equation}
The basis functions may be nonlinear transformations of the aggregate state, but the approximation is linear in its coefficients. If a home option $o=(h_k,w)$ is accepted, then $n_k$ and $n_k^H$ increase by one, the load of the corresponding home cell--slot pair increases or a new pair is activated, and $C_k$ is updated through the insertion operator. If a locker option $o=\ell$ is accepted, then $n_k$ increases by one, $n_k^H$ is unchanged, the load of locker $\ell$ increases or a new locker is activated, and $C_k$ is updated. If the customer opts out, only the epoch advances. Consequently, accepted-versus-opt-out value differences are computed by evaluating how a feasible option changes the aggregate booking state carried into the future.

\subsection{Offline learning of the post-decision value}
\label{subsec:vfa-learning}

The value model is trained offline using simulated booking trajectories. Let $\pi^R$ denote the rollout policy used to generate continuation-value labels. In the main implementation, $\pi^R$ is the calibrated fine-static policy because it is stable, competitive, and samples states in the region of the state space that the anchored dynamic policy is expected to visit. The seed sets used for static calibration, value-function training, validation, and final evaluation are disjoint.

For each training seed, we generate a complete arrival stream and the corresponding customer-choice random numbers. We simulate the booking horizon under $\pi^R$ and sample post-decision booking states along the trajectory. For a sampled booking state $B_i$ at epoch $k_i$, we clone the state and complete the remaining horizon under $\pi^R$. Let $\widetilde B^i_{\tau_i}$ be the terminal booking state reached by this rollout completion, where $\tau_i$ is the realized terminal epoch. Let $\operatorname{Rev}(B)$ denote the cumulative accepted-order revenue earned up to booking state $B$, and let $N_i^{\mathrm{out}}$ denote the number of opt-out outcomes observed during the cloned rollout completion after $B_i$. The continuation-profit label is
\begin{equation}
Y_i =
\underbrace{\left[\operatorname{Rev}\!\left(\widetilde B^i_{\tau_i}\right)-\operatorname{Rev}(B_i)\right]}_{\text{accepted-order revenue earned after }B_i}
-
\underbrace{\pi_{\mathrm{out}}N_i^{\mathrm{out}}}_{\text{opt-out penalties after }B_i}
-
\underbrace{C^{\mathrm{ful}}\!\left(\widetilde B^i_{\tau_i}\right)}_{\text{terminal fulfillment cost}}.
\label{eq:rollout-label}
\end{equation}

Equivalently, $Y_i$ is one Monte Carlo realization of the post-decision continuation value from $B_i$ under the rollout policy $\pi^R$. The label is future-only. The subtraction of $\operatorname{Rev}(B_i)$ removes all accepted-order revenue earned before the sampled state, so past revenue is not counted again. The term $N_i^{\mathrm{out}}$ counts only opt-outs that occur during the cloned continuation after $B_i$, so opt-out penalties incurred before the sampled state are not counted again. In contrast, the terminal cost is not differenced against a cost at $B_i$, because no fulfillment cost has yet been incurred at the post-decision state; the already accepted orders in $B_i$ still have to be served at the terminal stage. Thus, the learned value approximation represents future accepted-order revenue minus future opt-out penalties and terminal fulfillment cost, conditional on the sampled post-decision booking state. The labeled training set is $\mathcal{D}=\{(\phi(B_i),Y_i): i=1,\ldots,N\}$. The coefficients in \eqref{eq:vfa} are estimated from this labeled data using a linear regression estimator on the standardized basis. Appendix~\ref{app:offline_training} summarizes the offline training process.

\subsection{Candidate-menu construction and online SVAP--ADP policy}
\label{subsec:online-policy}
The exact action space in \eqref{eq:pricing-exact-opp} allows the provider to choose any insertion-feasible displayed menu. Enumerating all subsets of feasible AHD slots and locker options is computationally expensive and may also generate menus that are operationally irrelevant. We therefore implement SVAP--ADP on a restricted, state-dependent candidate-menu family.

Let $O_k^H$ be the feasible AHD options for the current customer. The displayed AHD sub-menu \(U_k^H\) contains the \(K^H\) feasible home slots with the smallest raw insertion costs. For OOH delivery, SVAP--ADP constructs a state-dependent locker candidate pool. The pool is built from three layers: the \(K^{\mathrm{near}}\) nearest accessible lockers, up to \(K^{\mathrm{act}}\) feasible active lockers already used in the current tentative plan, and up to \(K^{\mathrm{inact}}\) feasible inactive lockers with high forecast residual catchment demand. Duplicates are removed, and the total pool is capped at \(K^{\mathrm{pool}}\). SCAP-ADP further ranks this pool lexicographically by four keys: The first key places active lockers before inactive ones, the second favors low adjusted opportunity cost, the third breaks remaining ties toward proximity to the customer, and the fourth favors dense locker catchments. The nearest accessible locker is always retained. The displayed locker sub-menu \(U_k^L\) contains the top \(K^{\mathrm{offer}}\) lockers after this ranking. The complete displayed menu is $U_k = U_k^H \cup U_k^L$. Section~\ref{subsec:value_assortment_control} then varies the locker-pool design to isolate the additional value of OOH assortment control.

At each arrival, SVAP--ADP solves the approximate one-step problem over the feasible options and candidate menus generated by the construction procedure described above. It therefore prices only options that have already passed the insertion-feasibility filter and associated construction rules. The candidate-menu step restricts the assortment to this manageable feasible set. For each feasible option \(o\), the policy compares the opt-out successor \(B_{k+1}^{0}\) with the accepted successor \(B_{k+1}^{o}\), both constructed on cloned booking states. The estimated opportunity cost is
\begin{equation}
\widehat{\Delta}_k(S_k,o)
=
\widehat V^x(B_{k+1}^{0})
-
\widehat V^x(B_{k+1}^{o}).
\label{eq:svap-delta}
\end{equation}
Because \(\widehat V^x\) is a future-only post-decision value approximation, immediate revenue is not included in \eqref{eq:svap-delta}; it enters only the pricing objective below. To stabilize the learned value difference, SVAP--ADP blends it with the local insertion-cost signal:
\begin{equation}
\widetilde c_{\eta,k}(o)
=
\Pi_{[\underline c,\overline c]}
\left(
c_k^{\mathrm{ins}}(o)
+
\eta\left[
\widehat{\Delta}_k(S_k,o)-c_k^{\mathrm{ins}}(o)
\right]
\right),
\label{eq:adjusted-cost}
\end{equation}
where \(\eta\in[0,1]\) controls the ADP correction and \(\Pi_{[\underline c,\overline c]}\) denotes projection onto the clipping interval.

Prices are restricted to a trust-region grid around the fine-static anchor:
\begin{equation}
\mathcal G_k^{\mathrm{TR}}(o)
=
\left\{
p\in\mathcal G_{\kappa(o)}:
|p-p^F(o)|\le \tau_{\kappa(o)}^{\mathrm{TR}}
\right\}.
\label{eq:trust-region-grid}
\end{equation}
For each candidate menu \(U\in\mathcal M_k(S_k)\), SVAP--ADP solves
\begin{equation}
J_k(U)
=
\max_{p\in\prod_{o\in U}\mathcal G_k^{\mathrm{TR}}(o)}
\left\{
\sum_{o\in U}
q_k(o\mid S_k,(U,p))
\left[
r+p(o)-\widetilde c_{\eta,k}(o)
\right]
-q_k(0\mid S_k,(U,p))\pi_{out}
\right\},
\label{eq:anchored-pricing}
\end{equation}
and then deploys the menu-price pair with the largest score. The finite grids make \eqref{eq:anchored-pricing} small enough to solve by enumeration.

\begin{algorithm}[t]
\caption{Online SVAP--ADP decision at epoch \(k\)}
\label{alg:SVAP--ADP}
\begin{algorithmic}[1]
\State Construct feasible options \(\mathcal O_k(S_k)\), candidate menus \(\mathcal M_k(S_k)\), and the opt-out successor \(B_{k+1}^{0}\).
\For{\(o\in\mathcal O_k(S_k)\)}
    \State Construct the accepted successor \(B_{k+1}^{o}\) on a cloned tentative plan.
    \State Compute \(\widehat{\Delta}_k(S_k,o)\) by \eqref{eq:svap-delta} and \(\widetilde c_{\eta,k}(o)\) by \eqref{eq:adjusted-cost}.
\EndFor
\For{\(U\in\mathcal M_k(S_k)\)}
    \State Construct \(\mathcal G_k^{\mathrm{TR}}(o)\) for all \(o\in U\) and solve \eqref{eq:anchored-pricing} to obtain \(J_k(U)\) and \(p^*(U)\).
\EndFor
\State Select \(U_k^{\mathrm{SVAP}}\in\arg\max_{U\in\mathcal M_k(S_k)}J_k(U)\) and set \(p_k^{\mathrm{SVAP}}=p^*(U_k^{\mathrm{SVAP}})\).
\State Offer \((U_k^{\mathrm{SVAP}},p_k^{\mathrm{SVAP}})\), observe \(u_k\), and update the actual booking state to \(B_{k+1}^{u_k}\).
\end{algorithmic}
\end{algorithm}

Only the final update changes the actual system state; all earlier successor states are counterfactual evaluations. Thus, SVAP--ADP differs from a myopic insertion-cost policy by replacing \(c_k^{\mathrm{ins}}(o)\) with the clipped post-decision value signal \(\widetilde c_{\eta,k}(o)\), while the trust-region anchor keeps customer-facing prices close to the calibrated static benchmark. Algorithm~\ref{alg:SVAP--ADP} summarizes this online implementation.

It is worth mentioning that the online SVAP--ADP decision is a regularized approximation of the exact one-step opportunity-cost pricing problem. Three restrictions separate the implemented action from the full-information one-step optimizer. First, the candidate menu family $\mathcal M_k(S_k)$ restricts the assortment space. Second, the trust-region grid $G_k^{\mathrm{TR}}(o)$ restricts price deviations from the fine-static anchor. Third, the exact opportunity cost $\Delta_k(S_k,o)$ is replaced by the blended and clipped estimate $\tilde c_{\eta,k}(o)$ obtained from the learned post-decision value model. 

\section{Seattle Case Study}
\label{sec:seattle}
In this section, we evaluate SVAP--ADP on a real-world case based on the greater Seattle delivery network. Section~\ref{sec:benchmarks} introduces the benchmark policies. Section~\ref{sec:seattle_setup} describes the instance design. Sections~\ref{sec:value_dp}--\ref{sec:sense} report the experimental results.

\subsection{Benchmark Policies}
\label{sec:benchmarks}

We compare SVAP--ADP with five benchmarks designed to isolate the contribution of each modeling component: current-practice uniform pricing, static price calibration, myopic routing awareness, compact state aggregation, and anchored trust-region pricing. The first benchmark, Static-uniform, represents a common real-world pricing practice in which the provider posts one fixed home-delivery price and one fixed locker-delivery price throughout the booking horizon. The second benchmark, Static-fine, is the calibrated static anchor $\pi^F$ used by the anchored dynamic policies. All policies are evaluated on the same customer streams and choice-uniform streams, i.e., common uniform random-number streams used to sample realized MNL choices. Unless otherwise stated, they use the same feasible-option construction, locker candidate pool, and customer-choice model. Thus, performance differences are attributable to the pricing and opportunity-cost logic rather than to differences in the available assortment. Table~\ref{tab:svap_benchmark_design} summarizes the design.

\begin{table}[!htbp]
\centering
\caption{Benchmark design. Each policy is defined by its opportunity-cost estimate, value-function form, and pricing feasible set.}
\label{tab:svap_benchmark_design}
\footnotesize
\setlength{\tabcolsep}{4pt}
\renewcommand{\arraystretch}{1.02}
\begin{tabular*}{\textwidth}{@{\extracolsep{\fill}}lllll@{}}
\toprule
Policy & Opportunity cost & Value function & Pricing set & Isolates \\
\midrule
Static-uniform ($\pi^{\mathrm U}$) & State-independent & --- & Fixed $(p_h^{\mathrm U},p_\ell^{\mathrm U})$ & Current practice \\
Static-fine ($\pi^F$) & State-independent & --- & Fixed $p^F(o)$ & Static calibration \\
Myopic & Raw insertion cost & --- & Trust region $G_k^{\mathrm{TR}}(o)$ & Anticipation value \\
SVAP--ADP--6F & Learned VFA & 6-feature basis & Trust region $G_k^{\mathrm{TR}}(o)$ & State-basis richness \\
SVAP--ADP--FreePrice & Learned VFA & 17-component basis & Full grid $G_{\kappa(o)}$ & Anchoring value \\
SVAP--ADP & Learned VFA & 17-component basis & Trust region $G_k^{\mathrm{TR}}(o)$ & Full method \\
\bottomrule
\end{tabular*}
\end{table}

\begin{enumerate}
\item \emph{Static-uniform pricing ($\pi^{\mathrm U}$).} This policy represents a simple current-practice pricing rule used in many operational settings: the provider posts one fixed home-delivery price and one fixed locker-delivery price throughout the booking horizon. Formally, $p^{\mathrm U}(h,w)=p_h^{\mathrm U}$ for every home slot $w\in\mathcal S$ and $p^{\mathrm U}(\ell)=p_\ell^{\mathrm U}$ for every locker option $\ell$.

\item \emph{Static-fine pricing ($\pi^F$).} Slot-specific home surcharges and zone-segmented locker prices are calibrated offline on dedicated seeds, and the resulting price map $p^F(o)$ is applied at every arrival independently of the accumulated state. This is the fine-static anchor of SVAP--ADP; it absorbs the systematic price heterogeneity available to a policy that cannot observe the evolving consolidation and routing structure. Comparing Static-fine with Static-uniform quantifies the value of offline price segmentation and calibration, while comparing SVAP--ADP with Static-fine quantifies the additional value of dynamic state-dependent pricing.

\item \emph{Myopic pricing.} The opportunity-cost estimate collapses to the immediate raw insertion cost. Prices are optimized using the same anchored trust-region pricing structure as SVAP--ADP, but the value-function component is removed. Pricing therefore responds only to the current routing state, with no anticipation of future consolidation value or capacity scarcity. This policy isolates the contribution of forward-looking opportunity-cost estimation.

\item \emph{SVAP--ADP--6F.} The anticipation channel is active, but the value function is restricted to a compact six-feature aggregate representation defined in Eq.~\eqref{eq:core-aggregate-state}. Compared with SVAP--ADP, this benchmark changes only the value basis: it replaces the final 17-feature aggregate basis with the compact six-feature basis. All other components are kept identical to SVAP--ADP. This policy isolates the value of the richer state representation used by SVAP--ADP.

\item \emph{SVAP--ADP--FreePrice.} The final state-aggregation value function and the blended op\-por\-tu\-ni\-ty-cost estimate of SVAP--ADP are retained, but the trust-region anchor is removed: prices are optimized over the full mode-specific grids $G_{\kappa(o)}$ rather than the local grids $G_k^{\mathrm{TR}}(o)$. The learned correction is then free to post any feasible surcharge or discount, with no restriction to a neighborhood of the static anchor $p^F(o)$. With the opportunity-cost estimate held identical to SVAP--ADP, the single change is the pricing feasible set; this policy isolates the stabilizing contribution of the anchored trust region and exposes the customer-facing prices to value-approximation error.

\item \emph{SVAP--ADP.} The full method described in Section~\ref{sec:method}, combining the dimension-enriched aggregate state-value approximation, the blended opportunity-cost estimate, and anchored trust-region pricing.
\end{enumerate}

All policy comparisons use a paired Monte Carlo design. For each evaluation replication, we generate one customer-arrival stream and one stream of customer-choice random numbers, and then replay the same streams for every policy. Thus, each policy faces the same realized demand and latent choice shocks; only the offered menus, prices, and opportunity-cost logic differ. This common-random-number design reduces simulation noise and makes policy lifts interpretable as paired differences rather than as differences between unrelated simulations. The final performance evaluation uses 400 held-out Monte Carlo episodes. For policies with learned value functions, reported averages are taken over the same 400 episodes and 10 independently trained value models. We report mean episode outcomes and 95\% paired or nested confidence intervals for profit lifts.

The simulator is implemented in Python. Online insertion costs are produced by a cheapest-insertion heuristic; the offline routing backend adopts a Clarke--Wright savings method. In the computational experiments, $\widehat{\theta}$ is estimated by ridge regression to stabilize estimation under correlated basis terms. Results are reported as averages over 400 independent simulation episodes on a computer with a 1.9 GHz Intel Core i7-1370P processor and 32 GB of memory.

\subsection{Instance design}
\label{sec:seattle_setup}

\begin{figure}[htbp]
  \centering
  \includegraphics[width=0.92\textwidth]{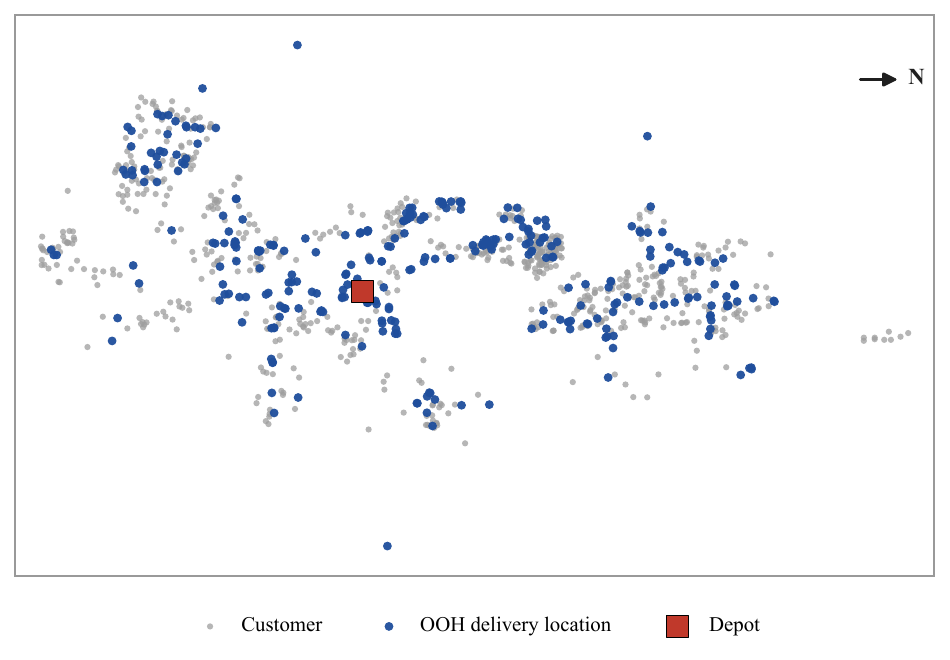}
  \caption{Seattle map}
  \label{fig:Seattle_map}
\end{figure}

The Seattle case is built from the publicly available Amazon dataset of \citet{merchan20242021}, which covers the greater Seattle metropolitan area and is geocoded by latitude and longitude. The network contains a single depot in the southwestern part of the metropolitan area near the main logistics corridor, 700 candidate customer home locations, and 299 candidate locker sites, as shown in Figure~\ref{fig:Seattle_map}. Each evaluation episode contains a deterministic booking horizon of $D=700$ sequential arrivals and $|\mathcal{S}|=3$ home-delivery time slots, corresponding to morning, afternoon, and evening. We set the vehicle-capacity parameter to $Q=120$, but hard fleet-capacity constraints are not enforced in the main Seattle setting. Accordingly, in the Seattle implementation, home-option feasibility is enforced through the time-window insertion check used by all policies, while fleet-capacity pressure is modeled economically through the route-count and fixed-cost terms in $C^{\mathrm{ful}}$ rather than through hard rejection. This common feasibility protocol is applied to every benchmark policy. 

Customer choices are generated from the MNL utility model~\eqref{eq:utility} with calibrated price sensitivity $\beta^p$ and locker-distance sensitivity $\beta^d$. Accepted orders earn base revenue $r$ plus the posted delivery-price adjustment, where positive adjustments are surcharges and negative adjustments are discounts. Realized fulfillment cost is computed at the end of each episode using a Clarke--Wright savings algorithm. In the main comparison, all policies use the same Seattle simulator, price grids, feasible-option construction, and candidate-menu protocol: each arriving customer is shown the two cheapest feasible home slots and up to six locker options from the shared ranked locker-candidate pool.

SVAP--ADP estimates aggregate continuation value and converts one-step state-value differences into opportunity costs, while stabilizing posted prices inside an anchored trust region around the fine-static policy $\pi^{\mathrm F}$. The full set of calibrated values, including choice-model coefficients, price grids, fulfillment-cost components, candidate-menu parameters, anchor prices, and key SVAP--ADP hyperparameters, is reported in Appendix~\ref{app:parameters}.

\subsection{Value of dynamic pricing}
\label{sec:value_dp}

We report results from 400 independent simulation episodes drawn from the held-out evaluation seed group. For policies with learned value functions, reported averages are taken over the same $400$ episodes and $10$ independently trained value models. Table~\ref{tab:seattle_profit_cost} reports profit, revenue, fulfillment cost, and demand split; Table~\ref{tab:seattle_mechanism} reports mechanism-level metrics and posted-price summaries. Profit lifts are reported relative to Static-uniform, a current-practice baseline that posts one fixed home-delivery price and one fixed locker-delivery price throughout the booking horizon. All policies use the same customer streams, choice-uniform streams, feasible-option construction, and locker candidate pool.

\begin{table}[htbp]
\centering
\caption{Profit, fulfillment cost, and demand split across policies on the Seattle case (mean over 400 episodes; lift relative to Static-uniform with 95\% paired or nested confidence interval).}
\label{tab:seattle_profit_cost}
\footnotesize
\setlength{\tabcolsep}{4pt}
\renewcommand{\arraystretch}{1.02}
\begin{tabular*}{\textwidth}{@{\extracolsep{\fill}}l l r r r r r r@{}}
\toprule
& \multicolumn{2}{c}{Profit} & & & \multicolumn{3}{c}{Demand split (\%)} \\
\cmidrule(lr){2-3}\cmidrule(lr){6-8}
Policy & Mean & \multicolumn{1}{c}{Lift (95\% CI)} & Revenue & Fulfillment & Home & Locker & Opt-out \\
\midrule
Static-uniform $(\pi^{\mathrm U})$& 6130.3 & --- & 15114.1 & 8792.5 & 21.4 & 73.2 & 5.5 \\
Static-fine $(\pi^{\mathrm F})$ & 6247.7 & $+117.4\;[105.6,\,129.2]$ & 15015.2 & 8578.3 & 21.2 & 73.4 & 5.4 \\
Myopic & 5716.9 & $-413.4\;[-435.7,\,-391.2]$ & 16456.6 & 10496.5 & 31.0 & 62.1 & 6.9 \\
SVAP--ADP-6F & 6022.9 & $-107.4\;[-179.8,\,-36.5]$ & 15865.0 & 9626.3 & 27.4 & 66.4 & 6.2 \\
SVAP--ADP-FreePrice & 5188.6 & $-941.7\;[-1420.7,\,-437.0]$ & 16677.3 & 11132.7 & 40.5 & 49.3 & 10.2 \\
\textbf{SVAP--ADP} & \textbf{6558.7} & $\mathbf{+428.4\;[413.0,\,443.2]}$ & \textbf{14987.1} & \textbf{8239.3} & \textbf{21.4} & \textbf{73.2} & \textbf{5.4} \\
\bottomrule
\end{tabular*}
\end{table}

\begin{table}[ht]
\centering
\caption{Operational mechanisms and posted prices across policies on the Seattle case (mean over 400 episodes).}
\label{tab:seattle_mechanism}
\small
\setlength{\tabcolsep}{6pt}
\renewcommand{\arraystretch}{1.05}
\begin{tabular*}{\textwidth}{@{\extracolsep{\fill}}lrrrrrr@{}}
\toprule
& \multicolumn{4}{c}{Mechanism diagnostics} & \multicolumn{2}{c}{Avg.\ price} \\
\cmidrule(lr){2-5}\cmidrule(lr){6-7}
Policy & Active $\ell$ & Singleton $\ell$ & Orders/active & Reuse & Home & Locker \\
\midrule
Static-uniform $(\pi^{\mathrm U})$ & 107.1 & 14.0 & 4.79 & 0.791 & 16.00 & $-1.00$ \\
Static-fine $(\pi^{\mathrm F})$ & 102.0 & 14.2 & 5.04 & 0.801 & 16.00 & $-1.17$ \\
Myopic & 92.0 & 13.0 & 4.73 & 0.788 & 15.53 & $+0.15$ \\
SVAP--ADP-6F & 87.1 & 10.4 & 5.35 & 0.812 & 15.56 & $-0.55$ \\
SVAP--ADP-FreePrice & 33.9 & 1.3 & 10.56 & 0.903 & 15.95 & $-1.14$ \\
\textbf{SVAP--ADP} & \textbf{86.7} & \textbf{8.1} & \textbf{5.92} & \textbf{0.831} & \textbf{15.97} & $\mathbf{-1.27}$ \\
\bottomrule
\end{tabular*}
\end{table}

Table~\ref{tab:seattle_profit_cost} first shows the value of moving beyond current-practice uniform pricing. Static-uniform posts a single home price and a single locker price, here $16$ and $-1$, across the entire booking horizon. Static-fine improves mean profit from $6130.3$ to $6247.7$, a lift of $117.4$, by using slot-specific home prices and density-segmented locker prices. The demand split changes little, but fulfillment cost falls from $8792.5$ to $8578.3$. Table~\ref{tab:seattle_mechanism} shows why: Static-fine uses fewer active lockers than Static-uniform and increases orders per active locker from $4.79$ to $5.04$, indicating that even offline segmentation improves consolidation relative to a one-price rule.

SVAP--ADP delivers the largest improvement over the current-practice baseline. It raises mean profit to $6558.7$, a lift of $428.4$ over Static-uniform. The gain is primarily operational rather than revenue-driven: SVAP--ADP has slightly lower revenue than Static-uniform, $14987.1$ versus $15114.1$, but reduces fulfillment cost substantially, from $8792.5$ to $8239.3$. The demand split is almost unchanged, with home share $21.4\%$, locker share $73.2\%$, and opt-out share $5.4\%$. Thus, the dynamic policy does not improve profit by suppressing demand or increasing customer charges; it improves profit by preserving the current-practice demand regime while steering accepted orders toward lower-cost fulfillment patterns.

The mechanism diagnostics confirm this interpretation. Compared with Static-uniform, SVAP--ADP reduces active lockers from $107.1$ to $86.7$ and singleton lockers from $14.0$ to $8.1$. Orders per active locker increase from $4.79$ to $5.92$, and locker reuse rises from $0.791$ to $0.831$. These changes indicate that SVAP--ADP concentrates OOH demand on a smaller and more heavily reused locker footprint. Average prices remain close to the static levels, with home price $15.97$ and locker price $-1.27$, so the learned policy acts as a state-dependent correction rather than an aggressive repricing rule.

The remaining benchmarks isolate why the full design is needed. Myopic pricing reacts only to immediate insertion costs. It raises revenue to $16456.6$, but fulfillment cost increases to $10496.5$, and profit falls below Static-uniform. The demand mix shifts sharply toward home delivery, with home share rising to $31.0\%$ and locker share falling to $62.1\%$, showing that immediate routing information alone is not enough. SVAP--ADP-6F activates the anticipation channel and improves the consolidation metrics relative to Myopic, but its six-feature basis is too coarse for the Seattle instance and remains below Static-uniform. SVAP--ADP-FreePrice keeps the learned value model but removes the anchored trust region; unrestricted price movement destabilizes the demand mix, raises opt-out, and produces the lowest profit.

Overall, the results show a clear progression. Static-uniform represents current practice with one fixed price per delivery mode. Static-fine improves on this baseline through offline price segmentation. SVAP--ADP adds state-dependent opportunity-cost pricing and produces the strongest result, improving profit by lowering terminal fulfillment cost and strengthening locker consolidation while maintaining the desired home--locker--opt-out mix.

\subsection{Value of Assortment Control}
\label{subsec:value_assortment_control}

The benchmark comparison in Tables~\ref{tab:seattle_profit_cost}--\ref{tab:seattle_mechanism} holds the locker candidate pool fixed across policies. We now relax this shared-pool comparison to quantify the value of assortment control. The baseline is the same Static-uniform, with prices $p_h=16$ and $p_\ell=-1$. We compare it with a calibrated Static-fine nearest-six benchmark and several SVAP--ADP candidate-pool variants. 

We parameterize an SVAP--ADP locker candidate pool by $(K^{\mathrm{pool}},K^{\mathrm{act}},K^{\mathrm{inact}},K^{\mathrm{offer}})$: the total candidate-pool size, active-locker cap, inactive-locker cap, and offered-locker cap. The near-six variant restricts SVAP--ADP to the same nearest-six locker menu as the static baselines. The sparse pool uses $(8,4,4,6)$, the active-heavy pool uses $(12,10,2,6)$, the reference pool $(12,6,6,6)$ is the frozen candidate-pool design used in the main benchmark, and the enriched pool uses $(16,8,8,8)$. This design separates the value of state-dependent pricing on a restricted menu from the additional value of dynamically expanding the locker assortment.

\begin{table}[!htbp]
\centering
\caption{Performance under restricted static menus on the Seattle case.}
\label{tab:seattle_assortment}
\footnotesize
\setlength{\tabcolsep}{2.5pt}
\renewcommand{\arraystretch}{1.02}
\begin{tabular*}{\textwidth}{@{\extracolsep{\fill}}l c r r r r r r@{}}
\toprule
& \multicolumn{2}{c}{Profit} & & & \multicolumn{3}{c}{Demand split (\%)} \\
\cmidrule(lr){2-3}\cmidrule(lr){6-8}
Policy & Mean & \multicolumn{1}{c}{Lift (95\% CI)} & Revenue & Fulfillment & Home & Locker & Opt-out \\
\midrule
Static-uniform, near-6 & 3722.1 & --- & 15025.8 & 11121.8 & 20.3 & 74.5 & 5.2 \\
Static-fine, near-6  & 3913.7 & $+191.6\;[174.7,\,208.5]$ & 15435.3 & 11315.8 & 23.0 & 71.1 & 5.9 \\
SVAP--ADP, near-6 & 4255.7 & $+533.6\;[417.3,\,651.3]$ & 15340.5 & 10884.8 & 22.5 & 71.8 & 5.7 \\
SVAP--ADP, sparse pool & 5683.0 & $+1960.9\;[1951.5,\,1971.3]$ & 15043.5 & 9165.2 & 22.1 & 72.4 & 5.6 \\
SVAP--ADP, active-heavy pool & 6331.0 & $+2608.9\;[2582.8,\,2633.8]$ & 15147.1 & 8610.5 & 23.5 & 70.6 & 5.9 \\
SVAP--ADP, reference pool & 6558.7 & $+2836.6\;[2819.9,\,2851.4]$ & 14987.1 & 8239.3 & 21.4 & 73.2 & 5.4 \\
\textbf{SVAP--ADP, enriched pool} & \textbf{6706.8} & $\mathbf{+2984.8\;[2969.2,\,3000.7]}$ & 14934.0 & \textbf{8059.7} & 18.8 & 76.4 & 4.8 \\
\bottomrule
\end{tabular*}
\end{table}

\begin{table}[!htbp]
\centering
\caption{Operational diagnostics under restricted static menus on the Seattle case.}
\label{tab:seattle_assortment_mechanism}
\footnotesize
\setlength{\tabcolsep}{4pt}
\renewcommand{\arraystretch}{1.02}
\begin{tabular*}{\textwidth}{@{\extracolsep{\fill}}l r r r r r r@{}}
\toprule
& \multicolumn{4}{c}{Locker utilization} & \multicolumn{2}{c}{Avg.\ price} \\
\cmidrule(lr){2-5}\cmidrule(lr){6-7}
Policy & Active $\ell$ & Singleton $\ell$ & Orders / active & Reuse & Home & Locker \\
\midrule
Static-uniform, near-6 & 206.0 & 71.7 & 2.53 & 0.605 & 16.00 & $-1.00$ \\
Static-fine, near-6 & 190.5 & 71.0 & 2.62 & 0.617 & 16.00 & $-0.64$ \\
SVAP--ADP, near-6 & 177.7 & 53.0 & 2.84 & 0.646 & 15.98 & $-0.75$ \\
SVAP--ADP, sparse pool & 113.5 & 19.8 & 4.47 & 0.776 & 15.97 & $-1.27$ \\
\textbf{SVAP--ADP, active-heavy pool} & \textbf{82.3} & \textbf{7.8} & \textbf{6.02} & \textbf{0.833} & 15.95 & $-1.32$ \\
SVAP--ADP, reference pool & 86.7 & 8.1 & 5.92 & 0.831 & 15.97 & $-1.27$ \\
SVAP--ADP, enriched pool & 102.7 & 11.7 & 5.21 & 0.808 & 15.98 & $-0.94$ \\
\bottomrule
\end{tabular*}
\end{table}

Table~\ref{tab:seattle_assortment} shows that the current-practice Static-uniform policy performs poorly when it is restricted to the nearest-six locker menu. It achieves mean profit $3722.1$ with fulfillment cost $11121.8$. Although its locker share is high at $74.5\%$, Table~\ref{tab:seattle_assortment_mechanism} shows that these locker orders are highly fragmented: the policy activates $206.0$ lockers, creates $71.7$ singleton lockers, and obtains only $2.53$ orders per active locker. Thus, a simple one-price rule with a nearest-six menu sends many customers to lockers, but it does not concentrate demand on lockers that are operationally useful.

Static-fine, near-6 improves profit by $191.6$ over Static-uniform. This gain comes from offline price calibration rather than from solving the underlying assortment problem. Static-fine raises revenue from $15025.8$ to $15435.3$, but fulfillment cost also increases from $11121.8$ to $11315.8$. Its mechanism diagnostics improve only modestly, with active lockers falling from $206.0$ to $190.5$ and reuse increasing from $0.605$ to $0.617$. This confirms that static price segmentation helps relative to current practice, but the restricted nearest-six menu still produces a fragmented locker footprint.

The SVAP--ADP near-six row isolates the value of state-dependent pricing without richer assortment. It uses the same nearest-six menu restriction as the static baselines, but increases profit by $533.6$ relative to Static-uniform, with confidence interval $[417.3,\,651.3]$. It also lowers fulfillment cost to $10884.8$ and improves all consolidation diagnostics relative to both static nearest-six policies: active lockers fall to $177.7$, singleton lockers fall to $53.0$, orders per active locker increase to $2.84$, and reuse rises to $0.646$. Dynamic opportunity-cost pricing therefore creates value even when the menu is restricted, although the nearest-six assortment remains a binding limitation.

Richer candidate pools produce much larger gains. The sparse pool raises the lift to $1960.9$ over Static-uniform and reduces fulfillment cost to $9165.2$. Its mechanism diagnostics show a sharp consolidation improvement: active lockers fall to $113.5$, singleton lockers to $19.8$, orders per active locker rise to $4.47$, and reuse reaches $0.776$. Even a modest dynamic candidate pool allows SVAP--ADP to move beyond nearest-distance offering and select lockers that are more valuable for system-wide consolidation. The active-heavy and reference pools deliver further improvements. The active-heavy pool achieves the tightest locker-utilization profile, with only $82.3$ active lockers, $7.8$ singleton lockers, $6.02$ orders per active locker, and reuse of $0.833$. The reference pool used in the main benchmark attains higher profit, $6558.7$, and a lift of $2836.6$ over Static-uniform. It balances reuse of active lockers with access to inactive high-potential lockers, producing strong consolidation while maintaining a more favorable demand split.

The enriched pool is the best performer in this sensitivity test. It reaches mean profit $6706.8$ and a lift of $2984.8$ over Static-uniform, with confidence interval $[2969.2,\,3000.7]$. Compared with the reference pool, it activates more lockers and has lower orders per active locker, but it reduces fulfillment cost further to $8059.7$, increases locker share to $76.4\%$, and lowers opt-out to $4.8\%$. This indicates that assortment richness creates value not only by reusing active lockers, but also by exposing additional lockers that are attractive to customers and efficient for the final route structure.

Overall, the experiment shows that offering decisions are a central source of value in mixed AHD--OOH delivery. Static-uniform represents current practice with one fixed price per delivery mode and a simple nearest-six locker menu; Static-fine improves on this baseline through offline price calibration. SVAP--ADP adds state-dependent opportunity-cost pricing, raising profit by $14.3\%$ when restricted to the same nearest-six menu. However, the gains become much larger once the policy can choose from richer candidate pools: the enriched-pool variant improves profit by $80.2\%$ relative to Static-uniform and by an additional $57.6\%$ relative to SVAP--ADP on the restricted menu. This indicates that assortment control is even more critical than pure within-menu price adjustment. Dynamic prices can steer customers toward operationally attractive options, but the candidate pool determines whether such high-value locker options are available in the first place.

\subsection{Robustness Across Demand--Capacity Regimes}
\label{subsec:demand_capacity_regimes}

We further test the robustness of SVAP--ADP across demand--capacity regimes by varying the fleet size $V\in\{10,18,30\}$ at fixed demand $D=600$ and per-vehicle capacity $Q=30$. This yields a stress regime with total capacity $300$, a balanced regime with total capacity $540$, and a slack regime with total capacity $900$. Hard capacity with capacity-full rejection is enforced throughout. Static-uniform, which posts one fixed home price and one fixed locker price, is used as the baseline in each regime. Static-fine is calibrated separately for each regime on dedicated tuning seeds, while all learned SVAP variants are evaluated on the same final locked customer streams and choice-uniform streams. Tables~\ref{tab:seattle_regimes} and~\ref{tab:seattle_regimes_mechanism} report performance and operational diagnostics; profit lifts are measured relative to the regime-specific Static-uniform policy.

\begin{table}[!htbp]
\centering
\caption{Performance across demand--capacity regimes on the Seattle case.}
\label{tab:seattle_regimes}
\footnotesize
\setlength{\tabcolsep}{4pt}
\renewcommand{\arraystretch}{1.02}
\begin{tabular*}{\textwidth}{@{\extracolsep{\fill}}l c r r r r r r@{}}
\toprule
& \multicolumn{2}{c}{Profit} & & & \multicolumn{3}{c}{Demand split (\%)} \\
\cmidrule(lr){2-3}\cmidrule(lr){6-8}
Policy & Mean & \multicolumn{1}{c}{Lift (95\% CI)} & Revenue & Fulfillment & Home & Locker & Opt-out \\
\midrule
\multicolumn{8}{l}{\emph{Stress regime ($V=10$, total capacity 300)}} \\
Static-uniform & 1842.5 & --- & 6867.4 & 4937.3 & 11.4 & 38.6 & 2.9 \\
Static-fine ($\pi^{\mathrm F}$) & 1907.4 & $+64.8\;[55.5,\,74.2]$ & 6808.2 & 4815.2 & 11.3 & 38.7 & 2.9 \\
Myopic & 1906.1 & $+63.5\;[48.9,\,78.2]$ & 7408.1 & 5394.5 & 15.4 & 34.6 & 3.6 \\
SVAP--ADP--6F & 1921.5 & $+78.9\;[54.0,\,102.8]$ & 7159.3 & 5144.8 & 14.4 & 35.6 & 3.1 \\
SVAP--ADP--FreePrice & 2018.3 & $+175.7\;[41.5,\,308.0]$ & 7729.7 & 5560.9 & 20.1 & 29.9 & 5.0 \\
\textbf{SVAP--ADP} & \textbf{2121.3} & $\mathbf{+278.8\;[265.5,\,291.1]}$ & \textbf{6856.3} & \textbf{4646.6} & \textbf{11.7} & \textbf{38.3} & \textbf{2.9} \\
\addlinespace
\midrule
\multicolumn{8}{l}{\emph{Balanced regime ($V=18$, total capacity 540)}} \\
Static-uniform & 4749.7 & --- & 12339.7 & 7434.2 & 20.4 & 69.6 & 5.2 \\
Static-fine ($\pi^{\mathrm F}$) & 4847.1 & $+97.4\;[86.5,\,108.2]$ & 12246.7 & 7245.8 & 20.2 & 69.8 & 5.1 \\
Myopic & 4586.0 & $-163.7\;[-183.2,\,-144.1]$ & 13558.9 & 8774.1 & 29.3 & 60.7 & 6.6 \\
SVAP--ADP--6F & 4603.1 & $-146.6\;[-206.7,\,-93.1]$ & 13164.4 & 8384.4 & 27.5 & 62.5 & 5.9 \\
SVAP--ADP--FreePrice & 4553.9 & $-195.8\;[-541.6,\,144.7]$ & 13639.3 & 8822.2 & 35.0 & 53.9 & 8.8 \\
\textbf{SVAP--ADP} & \textbf{5088.8} & $\mathbf{+339.1\;[262.4,\,385.5]}$ & \textbf{12253.6} & \textbf{7009.8} & \textbf{20.5} & \textbf{69.5} & \textbf{5.2} \\
\addlinespace
\midrule
\multicolumn{8}{l}{\emph{Slack regime ($V=30$, total capacity 900)}} \\
Static-uniform & 5093.5 & --- & 12961.3 & 7704.3 & 21.4 & 73.2 & 5.4 \\
Static-fine ($\pi^{\mathrm F}$) & 5197.2 & $+103.8\;[92.5,\,115.0]$ & 12873.1 & 7514.2 & 21.2 & 73.4 & 5.4 \\
Myopic & 4802.0 & $-291.5\;[-312.6,\,-270.4]$ & 14046.1 & 9038.0 & 30.4 & 62.7 & 6.9 \\
SVAP--ADP--6F & 4936.1 & $-157.4\;[-226.1,\,-95.8]$ & 13671.9 & 8551.9 & 28.2 & 65.7 & 6.1 \\
SVAP--ADP--FreePrice & 4479.9 & $-613.5\;[-920.1,\,-323.8]$ & 14262.6 & 9481.8 & 39.9 & 50.0 & 10.0 \\
\textbf{SVAP--ADP} & \textbf{5480.2} & $\mathbf{+386.8\;[363.3,\,405.0]}$ & \textbf{12867.8} & \textbf{7224.7} & \textbf{21.5} & \textbf{73.1} & \textbf{5.4} \\
\bottomrule
\end{tabular*}
\end{table}

\begin{table}[!htbp]
\centering
\caption{Operational diagnostics across demand--capacity regimes on the Seattle case.}
\label{tab:seattle_regimes_mechanism}
\footnotesize
\setlength{\tabcolsep}{4pt}
\renewcommand{\arraystretch}{1.02}
\begin{tabular*}{\textwidth}{@{\extracolsep{\fill}}l r r c r r r@{}}
\toprule
& \multicolumn{4}{c}{Locker utilization} & \multicolumn{2}{c}{Avg.\ price} \\
\cmidrule(lr){2-5}\cmidrule(lr){6-7}
Policy & Active $\ell$ & Singleton $\ell$ & Orders / active & Reuse & Home & Locker \\
\midrule
\multicolumn{7}{l}{\emph{Stress regime}} \\
Static-uniform & 81.3 & 22.1 & 2.85 & 0.648 & 16.00 & $-1.00$ \\
Static-fine ($\pi^{\mathrm F}$) & 78.3 & 21.6 & 2.97 & 0.663 & 16.00 & $-1.18$ \\
Myopic & 69.2 & 18.0 & 3.00 & 0.666 & 15.64 & $-0.18$ \\
SVAP--ADP--6F & 65.8 & 15.2 & 3.25 & 0.692 & 15.42 & $-0.81$ \\
SVAP--ADP--FreePrice & 40.6 & 5.8 & 5.12 & 0.788 & 15.93 & $-0.92$ \\
\textbf{SVAP--ADP} & \textbf{67.8} & \textbf{14.5} & \textbf{3.40} & \textbf{0.705} & \textbf{15.97} & $\mathbf{-1.16}$ \\
\addlinespace
\midrule
\multicolumn{7}{l}{\emph{Balanced regime}} \\
Static-uniform & 100.3 & 15.8 & 4.17 & 0.760 & 16.00 & $-1.00$ \\
Static-fine ($\pi^{\mathrm F}$) & 95.7 & 15.9 & 4.38 & 0.772 & 16.00 & $-1.17$ \\
Myopic & 86.4 & 14.1 & 4.22 & 0.763 & 15.57 & $0.06$ \\
SVAP--ADP--6F & 82.3 & 12.4 & 4.56 & 0.780 & 15.41 & $-0.48$ \\
SVAP--ADP--FreePrice & 45.5 & 3.4 & 8.20 & 0.867 & 15.95 & $-1.13$ \\
\textbf{SVAP--ADP} & \textbf{83.9} & \textbf{10.4} & \textbf{5.00} & \textbf{0.799} & \textbf{15.98} & $\mathbf{-1.22}$ \\
\addlinespace
\midrule
\multicolumn{7}{l}{\emph{Slack regime}} \\
Static-uniform & 102.0 & 15.4 & 4.31 & 0.768 & 16.00 & $-1.00$ \\
Static-fine ($\pi^{\mathrm F}$) & 97.3 & 15.6 & 4.53 & 0.779 & 16.00 & $-1.17$ \\
Myopic & 87.4 & 13.9 & 4.31 & 0.767 & 15.57 & $0.08$ \\
SVAP--ADP--6F & 83.3 & 11.8 & 4.74 & 0.789 & 15.46 & $-0.52$ \\
SVAP--ADP--FreePrice & 34.3 & 1.7 & 9.44 & 0.890 & 15.97 & $-1.13$ \\
\textbf{SVAP--ADP} & \textbf{83.3} & \textbf{9.2} & \textbf{5.28} & \textbf{0.810} & \textbf{15.98} & $\mathbf{-1.24}$ \\
\bottomrule
\end{tabular*}
\end{table}

Table~\ref{tab:seattle_regimes} shows that SVAP--ADP delivers a positive and statistically significant profit lift over Static-uniform in every regime. The lift is $278.8$ in the stress regime, $339.1$ in the balanced regime, and $386.8$ in the slack regime, with all confidence intervals bounded away from zero. Static-fine also improves over Static-uniform in all regimes, but the gains are smaller: $64.8$, $97.4$, and $103.8$, respectively. Thus, offline static calibration improves over a simple current-practice one-price rule, but the larger gains come from state-dependent opportunity-cost pricing.

The source of SVAP--ADP's improvement is primarily operational. In the stress regime, SVAP--ADP reduces fulfillment cost from $4937.3$ under Static-uniform to $4646.6$, while preserving a very similar demand split. In the balanced regime, fulfillment cost falls from $7434.2$ to $7009.8$, again with nearly unchanged home, locker, and opt-out shares. In the slack regime, fulfillment cost falls from $7704.3$ to $7224.7$. Across regimes, revenue remains close to the static baselines, so the profit gain is not driven by extracting more revenue or suppressing demand; it is driven by serving accepted demand more efficiently.

Table~\ref{tab:seattle_regimes_mechanism} confirms that the cost reduction is associated with stronger locker consolidation. Relative to Static-uniform, SVAP--ADP reduces singleton lockers from $22.1$ to $14.5$ in the stress regime, from $15.8$ to $10.4$ in the balanced regime, and from $15.4$ to $9.2$ in the slack regime. Orders per active locker increase from $2.85$ to $3.40$, from $4.17$ to $5.00$, and from $4.31$ to $5.28$, respectively. Locker reuse also rises in all regimes. The same consolidation mechanism therefore appears under tight, balanced, and slack capacity.

The myopic benchmark is not robust. Under stress, Myopic improves modestly over Static-uniform because capacity scarcity makes immediate insertion costs informative. In the balanced and slack regimes, however, Myopic falls below Static-uniform by $163.7$ and $291.5$, respectively. The failure mode is a demand shift toward home delivery and higher fulfillment cost. In the slack regime, for example, Myopic raises home share to $30.4\%$, lowers locker share to $62.7\%$, and increases fulfillment cost to $9038.0$. Immediate routing costs alone do not provide a sufficiently anticipatory signal for mixed AHD--OOH pricing.

The SVAP--ADP--6F benchmark captures part of the consolidation mechanism, especially in the stress regime, where it improves over Static-uniform by $78.9$. In the balanced and slack regimes, however, the compact six-feature basis is not sufficient: the policy loses $146.6$ and $157.4$ relative to Static-uniform. Although its mechanism diagnostics are generally better than Myopic, its demand mix remains too home-intensive compared with the full SVAP--ADP policy. This shows that the richer 17-component dimension-based representation is needed for stable opportunity-cost pricing across regimes.

The FreePrice ablation shows why the anchored trust region remains important. FreePrice can improve over Static-uniform in the stress regime, but its confidence interval is wide. In the balanced regime, its confidence interval includes zero, and in the slack regime it performs substantially worse than Static-uniform. The failure is caused by unstable demand shifts: home share rises to $35.0\%$ in the balanced regime and $39.9\%$ in the slack regime, while opt-out also increases. Although FreePrice reports high orders per active locker and high reuse, this reflects an overly narrow locker footprint combined with a large shift away from locker demand, not a desirable operating regime.

Overall, the regime experiment supports three conclusions. First, Static-uniform provides a current-practice baseline, and Static-fine improves on it through regime-specific static calibration. Second, SVAP--ADP improves substantially over Static-uniform in all demand--capacity regimes, with gains driven by lower fulfillment cost and stronger locker consolidation rather than by demand suppression. Third, the design choices of SVAP--ADP remain important across regimes: the richer value basis outperforms the compact six-feature approximation, and anchored trust-region pricing prevents the instability observed under unrestricted FreePrice optimization.

\subsection{Sensitivity analysis}
\label{sec:sense}

The previous experiments show that SVAP--ADP delivers consistent profit gains across benchmark policies, assortment designs, and demand--capacity regimes. We now examine the robustness of these gains to three key parameters. Section~\ref{subsec:gamma_sensitivity} varies the fulfillment multiplier $\gamma$, which controls the weight of routing cost in the objective. Section~\ref{subsec:eta_sensitivity} varies the ADP-correction scale $\eta$, which interpolates between immediate raw insertion costs and learned state-value opportunity costs. Section~\ref{subsec:trust_region_sensitivity} varies the home and locker trust-region radii that bound customer-facing price deviations around the static anchor.

\subsubsection{Sensitivity to the fulfillment multiplier}
\label{subsec:gamma_sensitivity}

\begin{figure}[!htbp]
  \centering
  \includegraphics[width=0.92\textwidth]{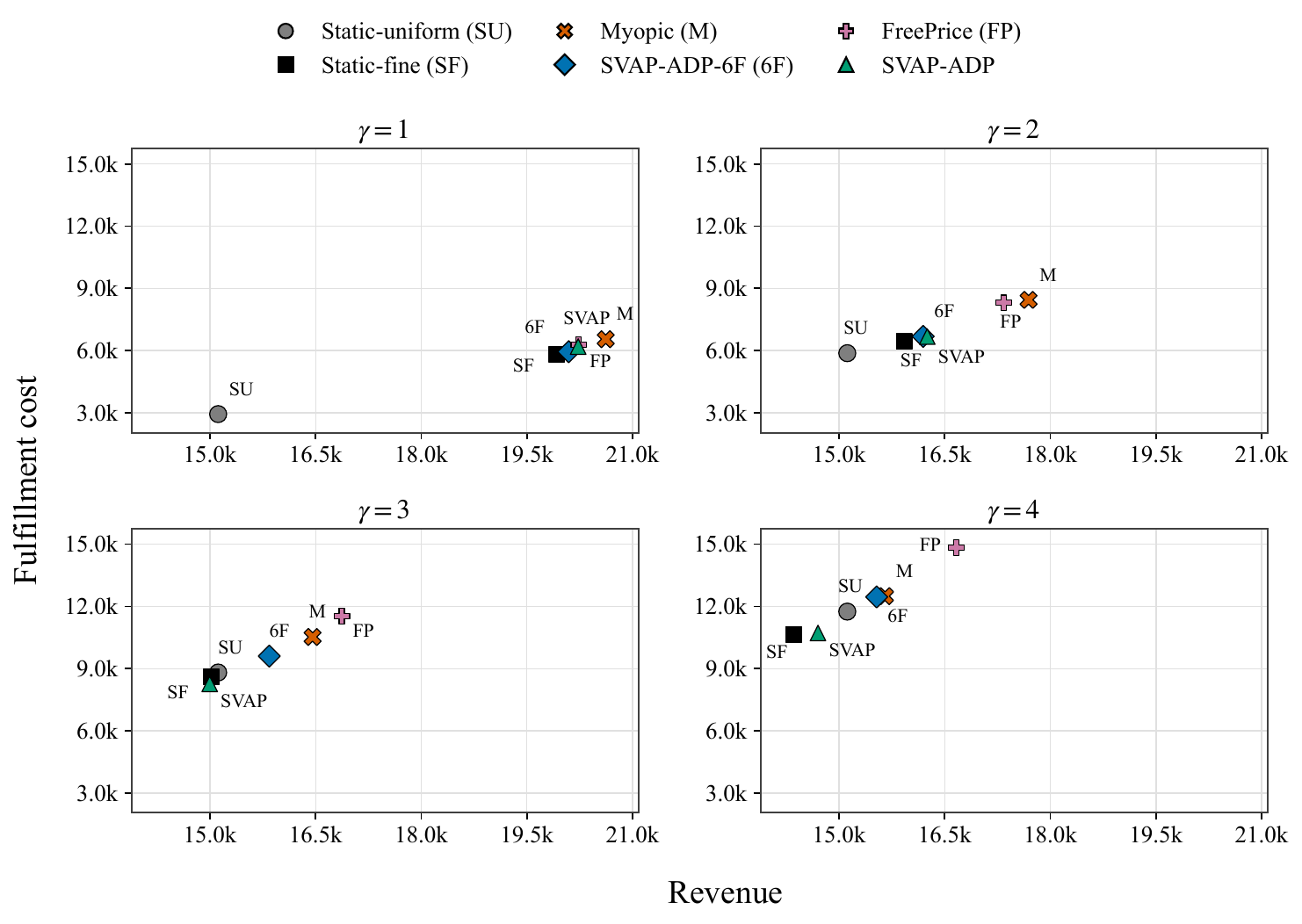}
  \caption{Revenue--fulfillment frontier across fulfillment-cost multipliers $\gamma\in\{1,2,3,4\}$.}
  \label{fig:rev_fulfill_frontier}
\end{figure}

\begin{figure}[htbp]
  \centering
  \includegraphics[width=0.92\textwidth]{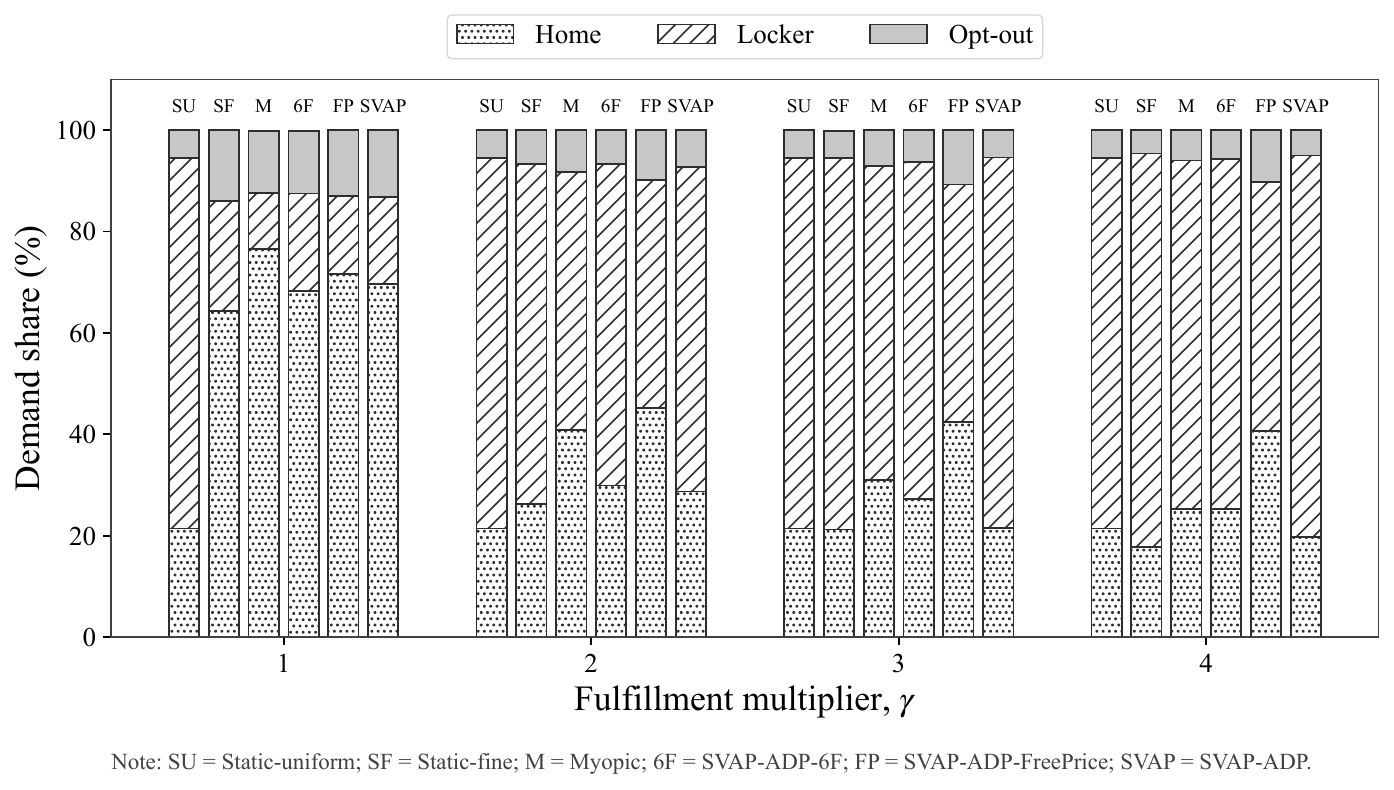}
  \caption{Demand-share decomposition across fulfillment-cost multipliers $\gamma\in\{1,2,3,4\}$. Bars report mean home, locker, and opt-out shares over $200$ final evaluation episodes.}
  \label{fig:demand_mix_gamma}
\end{figure}

We vary the fulfillment multiplier $\gamma\in\{1,2,3,4\}$ to scale the operational importance of routing cost relative to base revenue. The setting ranges from comparatively inexpensive fulfillment ($\gamma=1$) to a regime in which routing cost dominates the pricing decision ($\gamma=4$). This sensitivity analysis uses $200$ final locked evaluation episodes for each value of $\gamma$. Figures~\ref{fig:rev_fulfill_frontier} and~\ref{fig:demand_mix_gamma} report the revenue--fulfillment frontier and demand-mix decomposition. 

Figure~\ref{fig:rev_fulfill_frontier} shows that Static-uniform provides a useful reference point for how simple fixed pricing behaves as fulfillment costs change. Because SU posts the same prices at every $\gamma$, its revenue and demand mix are essentially unchanged, while its fulfillment cost scales mechanically with $\gamma$: profit falls from $11986.4$ at $\gamma=1$ to $3173.8$ at $\gamma=4$. Static-fine improves over SU at every multiplier by recalibrating static prices, but the relative position of SVAP--ADP shows that state-dependent opportunity-cost pricing becomes especially valuable once routing costs are economically important.

When $\gamma=1$, fulfillment is cheap and the policies cluster in the high-revenue region of the frontier. SU is much more conservative, with lower revenue and lower fulfillment cost, while Static-fine and the dynamic policies accept a more home-intensive, high-revenue demand mix. SVAP--ADP improves over SU by $1586.5$. However, the frontier also shows that when fulfillment is inexpensive, aggressive revenue extraction is not heavily penalized by routing cost; consequently, the advantage of the full SVAP--ADP structure over other dynamic variants is less pronounced in this regime.

At $\gamma=2$, the policies begin to separate along the revenue--fulfillment trade-off. SVAP--ADP remains close to the Static-fine region of the frontier but reduces unnecessary fulfillment burden relative to the more aggressive Myopic and FreePrice policies. Its lift over SU is $271.6$, with confidence interval $[167.0,\,369.0]$. This intermediate case shows that static calibration already captures part of the value of moving away from the fixed-price SU baseline, while the dynamic value signal starts to provide additional operational discipline.

The value of SVAP--ADP becomes clearer at $\gamma=3$ and $\gamma=4$. At $\gamma=3$, SVAP--ADP improves over SU by $432.4$, with confidence interval $[417.8,\,446.5]$. At $\gamma=4$, the lift increases to $633.2$, with confidence interval $[617.1,\,648.1]$. In both cases, SVAP--ADP is positioned near the lower-cost part of the frontier: it gives up some revenue relative to Myopic and FreePrice, but avoids their large fulfillment-cost penalties. This confirms that the learned opportunity-cost signal is most valuable when routing cost is sufficiently important for consolidation to determine profitability.

The frontier also highlights the failure modes of the ablations. Myopic pricing moves toward higher revenue but also much higher fulfillment cost, especially as $\gamma$ increases. FreePrice is even less stable: without the anchored trust region, it moves to the high-revenue, high-cost corner of the frontier and performs poorly when fulfillment is expensive. The six-feature variant is more disciplined than Myopic and FreePrice, but it does not match the full SVAP--ADP policy at high multipliers. Thus, the frontier supports the two main design choices of SVAP--ADP: a richer value representation and anchored trust-region pricing.

Figure~\ref{fig:demand_mix_gamma} explains the frontier through demand composition. Static-uniform keeps an almost fixed demand mix across all multipliers, with approximately $21.4\%$ home, $73.1\%$ locker, and $5.5\%$ opt-out. This is expected because SU does not adapt prices to the fulfillment multiplier. Static-fine and SVAP--ADP, by contrast, adjust the pricing logic to the cost environment. When $\gamma=1$, they allow a much larger home share because home delivery is relatively cheap. When $\gamma$ increases, both policies shift demand toward lockers and reduce home delivery.

At $\gamma=3$, SVAP--ADP has a demand mix close to SU, with home share $21.5\%$, locker share $73.1\%$, and opt-out $5.4\%$, but it achieves higher profit by improving the operational structure of the accepted orders. This indicates that the gain is not simply from changing the aggregate home--locker split; it comes from state-dependent pricing that steers demand toward lockers and home slots that are more favorable for final route construction. At $\gamma=4$, SVAP--ADP pushes further toward lockers, with home share $19.7\%$ and locker share $75.3\%$, which is appropriate when fulfillment cost is highest.

The contrast with Myopic and FreePrice is sharp. Myopic remains too home-oriented as fulfillment becomes expensive: at $\gamma=3$, its home share is $31.0\%$, compared with $21.5\%$ under SVAP--ADP, and at $\gamma=4$ it remains at $25.3\%$. FreePrice is more distorted. At $\gamma=4$, it produces a $40.7/49.1/10.3$ home--locker--opt-out split, moving demand in the wrong direction and roughly doubling opt-out relative to SVAP--ADP. This confirms that unrestricted price movement can destabilize customer-facing decisions when routing costs are large.

The six-feature variant provides an intermediate case. It captures part of the locker-steering mechanism, but its compact state representation is not sufficient at high multipliers. At $\gamma=4$, SVAP--ADP--6F has a home share of $25.3\%$ and locker share of $69.0\%$, whereas full SVAP--ADP has $19.7\%$ home and $75.3\%$ locker. The richer 17-component state basis is therefore needed to stabilize opportunity-cost pricing when the fulfillment-cost term dominates the objective.

Overall, the multiplier sensitivity supports three conclusions. First, Static-uniform is a useful current-practice baseline: it maintains a fixed demand regime, but it cannot adapt to changing fulfillment economics. Second, static calibration improves on SU at every multiplier, but the larger and more robust gains come from SVAP--ADP, especially when $\gamma\geq 3$. Third, the two main SVAP--ADP design elements remain important across the frontier: the richer state-value approximation prevents the compact policy from drifting toward costly home-heavy demand, and the anchored trust region prevents the FreePrice policy from overreacting through extreme demand shifts and higher opt-outs.

\subsubsection{Sensitivity to the ADP-correction scale}
\label{subsec:eta_sensitivity}

\begin{figure}[htbp]
  \centering
  \includegraphics[width=0.97\textwidth]{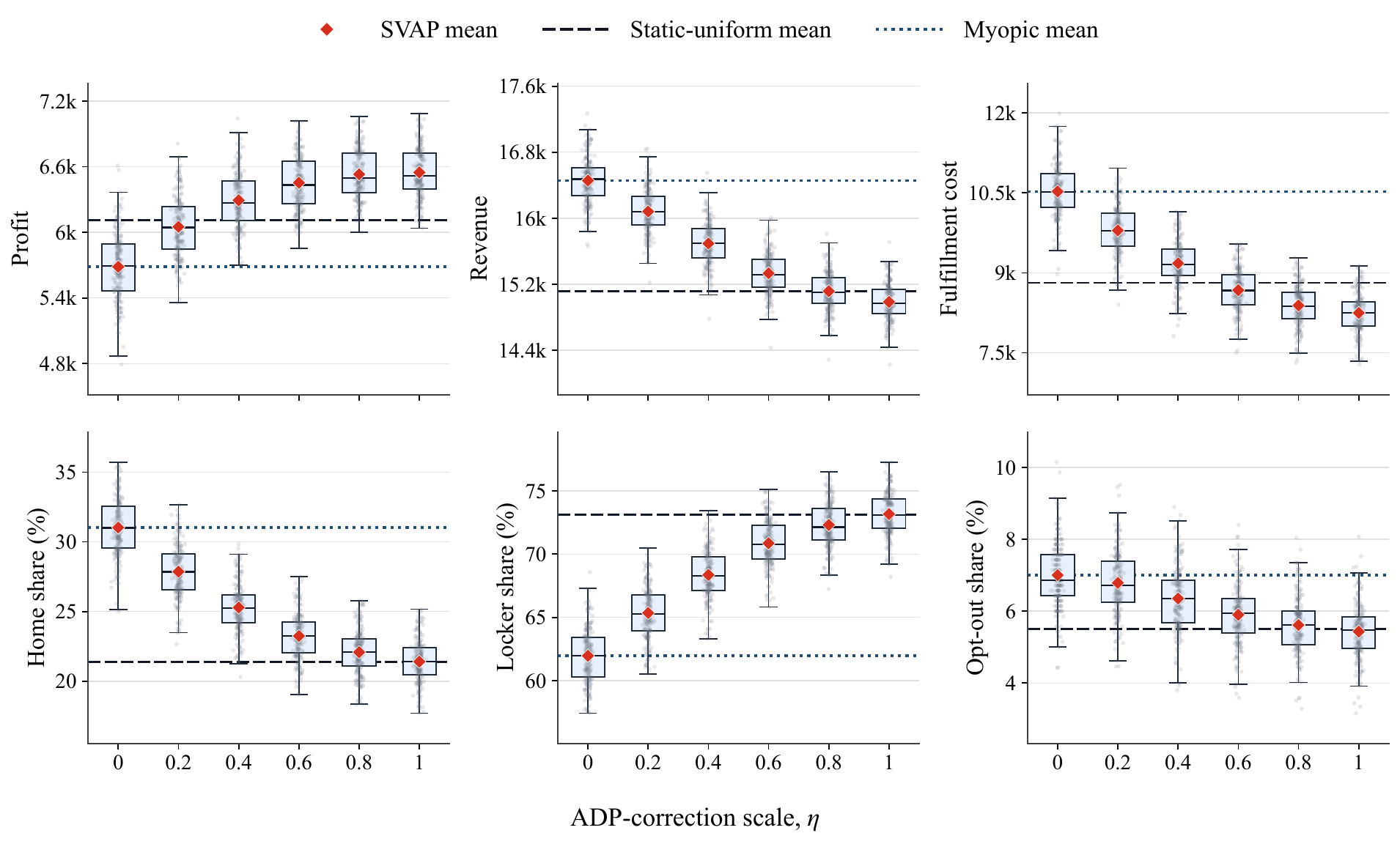}
  \caption{Sensitivity to the ADP-correction scale $\eta$. Red diamonds mark SVAP means. Dashed and dotted horizontal lines mark the Static-uniform and Myopic means for each metric, respectively.}
  \label{fig:eta_sensitivity}
\end{figure}

We next vary the ADP-correction scale $\eta$ in~\eqref{eq:adjusted-cost}. This parameter controls the interpolation between the immediate raw insertion-cost signal and the learned continuation-value opportunity cost. At $\eta=0$, the adjusted cost collapses to the raw insertion cost and the policy coincides with the myopic endpoint. At $\eta=1$, the adjusted cost fully uses the learned value-based correction before clipping. We evaluate $\eta\in\{0,0.2,0.4,0.6,0.8,1.0\}$ using the 200-seed final sensitivity protocol. Figure~\ref{fig:eta_sensitivity} reports episode-level distributions of profit, revenue, fulfillment cost, and demand split. The dashed horizontal reference is Static-uniform; the dotted reference is the myopic mean.

The profit panel shows a clear improvement as the ADP correction becomes stronger. At $\eta=0$, mean profit is $5685.5$, matching the myopic endpoint and lying well below the Static-uniform mean of $6111.3$. Increasing $\eta$ to $0.2$ raises profit to $6050.1$, nearly closing the gap to Static-uniform. At $\eta=0.4$, profit reaches $6293.0$ and moves above the current-practice baseline. Further increases continue to improve performance, with mean profits of $6453.1$, $6530.9$, and $6547.5$ at $\eta=0.6$, $0.8$, and $1.0$, respectively. The gains flatten near the top of the grid, indicating that the method is not sharply tuned to a single value of $\eta$; performance is strong over a broad range once the learned value signal receives sufficient weight.

The revenue and fulfillment-cost panels explain the profit pattern. As $\eta$ increases, revenue decreases smoothly from the myopic level of $16456.9$ to $14985.7$, while fulfillment cost falls much more sharply from $10526.3$ to $8248.1$. Relative to Static-uniform, the final $\eta=1$ policy gives up only a small amount of revenue, but it produces a large fulfillment-cost reduction. This is the intended role of the ADP correction: it prevents the online pricing problem from favoring immediate revenue when the current acceptance would create costly downstream routing commitments.

The demand-mix panels show that the cost reduction comes from a controlled shift toward OOH consolidation rather than from demand suppression. At $\eta=0$, the policy is home-intensive, with home share $31.0\%$, locker share $62.0\%$, and opt-out $7.0\%$. As $\eta$ increases, home share declines steadily, locker share rises, and opt-out falls. At $\eta=1$, the demand split is approximately $21.4\%$ home, $73.2\%$ locker, and $5.4\%$ opt-out, which is almost identical to the Static-uniform split of $21.4\%$, $73.1\%$, and $5.5\%$. Thus, the profit gain over Static-uniform is not achieved by changing aggregate demand volume or inducing additional opt-outs. Instead, the learned correction preserves the same broad demand regime while improving which options are encouraged in each state.

Two observations are important. First, the effect of $\eta$ is asymmetric: moving from $\eta=0$ to moderate values yields large fulfillment savings, whereas the marginal gain beyond $\eta=0.8$ is small. This suggests that raw insertion cost contains useful local information, but it is too myopic and too home-biased on its own. The value-based correction accounts for future locker consolidation and route-dispersion effects that are not visible in the immediate insertion cost. Second, the high-$\eta$ policies remain stable. Even at $\eta=1$, opt-out stays close to the Static-uniform level and the demand mix remains balanced. This contrasts with the instability observed in the FreePrice ablation, and confirms that the trust-region anchor and clipping layer stabilize the learned value signal.

Overall, the $\eta$ sensitivity supports the blended adjusted-cost formula. Increasing $\eta$ systematically lowers fulfillment cost, raises profit, increases locker share, and reduces opt-out relative to the myopic endpoint. More importantly, the final policy improves substantially over the Static-uniform current-practice baseline while maintaining nearly the same aggregate demand split. We therefore retain $\eta=1.0$ in the main experiments and interpret $\eta$ as a stabilizing interpolation parameter between immediate insertion-cost pricing and full state-value opportunity-cost pricing.

\subsubsection{Sensitivity to trust-region parameters}
\label{subsec:trust_region_sensitivity}

\begin{figure}[t]
    \centering
    \includegraphics[width=\textwidth]{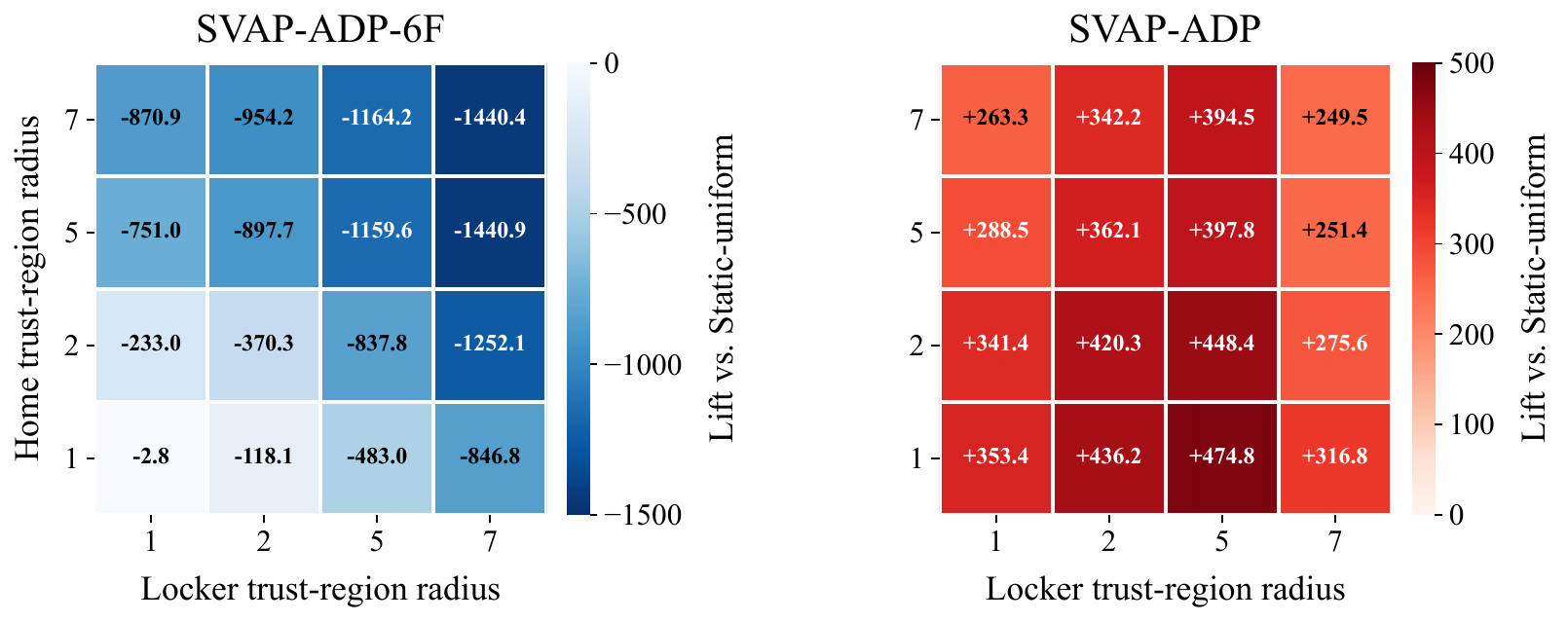}
    \caption{Profit lift over Static-uniform across home and locker trust-region grids $(\tau_h,\tau_\ell)\in\{1,2,5,7\}^2$.}
    \label{fig:trust_region_heatmap}
\end{figure}

We next assess the robustness of SVAP--ADP to the trust-region parameters that bound the per-option price deviation from the static anchor. Let $\tau_h$ and $\tau_\ell$ denote the home and locker trust-region radii. Larger values give the online pricing problem more freedom to depart from the anchor when converting opportunity-cost estimates into customer-facing prices. We vary each parameter in $\{1,2,5,7\}$ and report profit lift relative to Static-uniform. Because this is a sensitivity analysis, each cell is evaluated using the 200-seed sensitivity protocol. The resulting $4\times4$ heatmaps appear in Figure~\ref{fig:trust_region_heatmap}.

The left panel shows that SVAP--ADP--6F is highly sensitive to the trust-region widths. Relative to Static-uniform, all cells are non-positive, and the tightest setting $(\tau_h,\tau_\ell)=(1,1)$ is nearly tied with the baseline, with lift $-2.8$. Performance deteriorates quickly as either trust region widens, reaching losses below $-1{,}400$ in the widest locker-price settings. The decline is especially strong along the locker dimension: for any fixed home radius, increasing $\tau_\ell$ from $1$ to $7$ sharply reduces profit. This indicates that the compact six-feature value model does not provide a sufficiently stable opportunity-cost signal when the pricing layer is allowed to deviate substantially from the anchor.

The final SVAP--ADP policy behaves very differently. Every cell in the right panel produces a positive lift over Static-uniform, ranging from $249.5$ to $474.8$. The best cell is $(\tau_h,\tau_\ell)=(1,5)$, and the high-performing region is broad.  Moderate locker-price flexibility is consistently valuable.

The comparison between the two panels is the main message. The same trust-region grid that destabilizes the six-feature model remains safe for the 17-feature SVAP--ADP model. A coarse value representation can become harmful when the policy is given too much pricing freedom, because approximate opportunity costs are amplified into customer-facing price changes. By contrast, the richer 17-feature basis produces a more reliable state-value signal, allowing the policy to use price flexibility to exploit OOH consolidation opportunities while remaining profitable relative to the current-practice Static-uniform baseline.

The heatmap also confirms the role of anchoring. The final SVAP--ADP policy benefits from widening the locker trust region from $1$ to moderate values, but performance falls when $\tau_\ell=7$. This mirrors the FreePrice ablation: price flexibility is valuable only while the anchor still prevents excessive deviations. Overall, the sensitivity supports the final design of SVAP--ADP as an anchored value-based pricing policy. The trust region and the value basis are complementary: the trust region stabilizes customer-facing prices, and the richer value basis determines whether the policy can safely use that flexibility to improve profit.

\section{Conclusion}
\label{sec:conclusion}

This paper studied the DOPMDO problem, in which a logistics service provider dynamically selects and prices attended home-delivery time slots and out-of-home pickup options for sequentially arriving customers. We formulated the problem as a finite-horizon Markov decision process that links customer choice, immediate revenue, and route-dependent terminal fulfillment cost. To solve this problem, we proposed SVAP--ADP. The method approximates aggregate post-decision state value, converts accepted-versus-rejected state-value differences into option-level opportunity costs, and uses these opportunity costs to make online pricing decisions within an anchored trust region around a calibrated static policy. The main design principle is that dynamic pricing should exploit state-dependent consolidation opportunities, but should do so through stable, bounded price adjustments rather than unrestricted learned repricing.

The Seattle case study shows that SVAP--ADP generates robust profit gains over both current-practice and calibrated static benchmarks. Relative to Static-uniform, which posts one fixed home price and one fixed locker price throughout the booking horizon, SVAP--ADP increases mean profit by $428.4$ per episode, or approximately $7.0\%$; relative to the calibrated fine-static benchmark, the gain is $310.9$, or approximately $5.0\%$. The improvement is primarily operational: SVAP--ADP keeps revenue and the aggregate home--locker--opt-out split close to the static baselines, but reduces terminal fulfillment cost by improving locker consolidation and avoiding weakly used stops. The assortment-control experiments further show that offering decisions are a major source of value. When the menu is fixed to the nearest-six lockers, SVAP--ADP improves profit by $14.3\%$ relative to Static-uniform; when the candidate pool is enriched, the lift rises to $80.2\%$, an additional $57.6\%$ improvement over SVAP--ADP on the restricted menu. Thus, dynamic prices can only steer demand among the options made available, while the candidate pool determines whether the policy can expose lockers that are both attractive to customers and valuable for route consolidation.

The benchmark and ablation results clarify why the full SVAP--ADP design is needed. Static-fine improves over Static-uniform, showing the value of offline price segmentation, but SVAP--ADP adds state-dependent opportunity-cost pricing beyond static calibration. Myopic pricing raises revenue but creates costly downstream fulfillment patterns, the six-feature value-function benchmark is too coarse to produce stable opportunity-cost estimates, and the FreePrice ablation performs poorly because unrestricted price movements distort demand and increase fulfillment cost. These comparisons show that forward-looking value approximation, a sufficiently informative aggregate state basis, and anchored price stabilization are complementary components of the proposed method.

The sensitivity analyses show that the results are not driven by a narrow experimental setting. Across demand--capacity regimes, SVAP--ADP improves over the regime-specific current-practice baseline in the stress, balanced, and slack cases. The fulfillment-multiplier analysis shows that the value of the learned opportunity-cost signal increases when routing cost becomes economically important. The $\eta$ sensitivity shows that the blended adjusted-cost formula is not fragile: performance improves as the learned value signal receives more weight and then plateaus near the implementation value. The trust-region heatmaps further show that the final 17-feature SVAP--ADP model delivers positive lift throughout the tested grid, whereas the six-feature model remains weak or negative.

The managerial implication is that mixed AHD--OOH systems should not treat pricing, locker assortment, and routing as separate decisions. Simple fixed prices are easy to operate, and calibrated static prices improve on them, but both miss the state-dependent value of reusing active lockers, avoiding singleton pickup points, and managing route dispersion over the booking horizon. SVAP--ADP provides a practical way to incorporate these effects into online prices while remaining close to an interpretable static anchor. For logistics providers, this suggests that dynamic offering and pricing should be deployed as controlled corrections around reliable static tariffs, with candidate-pool design used to ensure that the policy has access to operationally meaningful consolidation opportunities.

Several limitations remain and suggest directions for future work. First, SVAP--ADP is an offline-trained, online-deployed policy.  This separation makes the online decision process fast and stable, but it also means that performance depends on the fidelity of the offline simulator, the routing surrogate, the choice model, and the representativeness of the training seeds. A natural extension is to update the value model and anchor periodically as new transaction and routing data become available. Second, the customer-choice model is calibrated rather than estimated from transaction-level booking data. The MNL specification provides a controlled basis for comparing policies, but deployment would require revealed-preference data that identify how customers trade off price, pickup distance, delivery mode, and time-window convenience. Future work could estimate richer choice models, such as nested-logit or latent-class models, to capture heterogeneous substitution patterns. Third, the current model treats each customer interaction as a one-shot decision. It does not capture repeated customer exposure to dynamic delivery prices, fairness perceptions, or long-run effects on loyalty and platform trust. Fourth, locker operations are simplified. The model focuses on delivery-side consolidation, but real locker systems also involve compartment availability, parcel-size compatibility, failed or delayed pickups, maintenance outages, and customer-service constraints at pickup locations. Extending SVAP--ADP to include these operational details would help assess whether routing-efficient locker recommendations remain feasible and attractive in real deployment.

\section*{Acknowledgments}
This work was supported by the European Commission and Swedish Energy Agency through the project E-Laas (Energy optimal urban Logistics As A Service).

\printbibliography
\clearpage
\appendix

\section{Notation and abbreviations}\label{app:notation}
Table~\ref{tab:key-notation-abbrev} summarizes the key notation and abbreviations used in the paper.

\begin{table}[H]
\centering
\caption{Key notation and abbreviations.}
\label{tab:key-notation-abbrev}
\scriptsize
\setlength{\tabcolsep}{2pt}
\renewcommand{\arraystretch}{0.90}
\begin{tabularx}{\textwidth}{@{}lY@{\hspace{1.25em}}lY@{}}
\toprule
\multicolumn{4}{@{}l}{\textit{Notation}}\\
\midrule
\multicolumn{2}{@{}l}{\textit{Problem setting and state}} & \multicolumn{2}{l}{\textit{Options, choice, and pricing}} \\
$\mathcal V=\{1,\ldots,V\}$ & Fleet of homogeneous vehicles & $O_k(S_k)$ & Feasible option set at decision state $S_k$ \\
$Q$ & Vehicle capacity & $\pi_{out}$ & Opt-out penalty \\
$G=(N^0,E)$ & Travel-time network & $(h_k,w)$, $w\in S$ & AHD option in time window $w$ \\
$L$, $S$ & Parcel-locker set; AHD time-window set & $\ell\in L_k$ & OOH locker option available to customer $k$ \\
$D$; $k$; $\tau$ & Maximum horizon length; decision epoch; terminal epoch & $U_k\subseteq O_k(S_k)$ & Displayed assortment/menu \\
$B_k=(k,\Omega_k,P_k)$ & Booking state after previous decisions & $x_k=(U_k,p_k,\{P_k^o\}_{o\in U_k})$ & Menu, prices, and option-contingent plans \\
$S_k=(B_k,C_k^{\mathrm{new}})$ & Decision state for the current customer & $p_k(o)$; $r$ & Price adjustment for option $o$; base revenue \\
$C_k^{\mathrm{new}}=(h_k,b_k,L_k)$ & Current customer information & $v_k(o)$ & Systematic utility of option $o$ \\
$h_k$, $b_k$, $L_k$ & Home location, load requirement, and candidate lockers & $\mu_w$, $\mu_\ell$ & AHD slot and OOH baseline utilities \\
$\Omega_k$, $P_k$ & Accepted-order set and tentative fulfillment plan & $\beta^p$, $\beta^d$ & Price and locker-distance sensitivities \\
$u_k\in U_k\cup\{0\}$ & Realized choice; $0$ is the outside option & $q_k(u\mid S_k,x_k)$ & MNL choice probability under decision $x_k$ \\
\addlinespace
\multicolumn{2}{@{}l}{\textit{Routing, costs, and value}} & \multicolumn{2}{l}{\textit{SVAP--ADP state and learning}} \\
$B_{k+1}^{0}$, $B_{k+1}^{o}$ & Opt-out and option-acceptance successor states & $X(B_k)$ & Core aggregate booking-state summary \\
$P_k^o$ & Tentative plan if option $o$ is accepted & $n_k$, $n_k^H$, $n_k^L$ & Total, AHD, and OOH accepted-order counts \\
$c_k^{\mathrm{ins}}(o)$ & Immediate marginal insertion cost & $A_k^H$, $A_k^L$ & Active AHD cell--slot map and active locker map \\
$C^{\mathrm{ful}}(B_\tau)$ & Terminal fulfillment cost & $m_{ak}^H$, $m_{\ell k}^L$ & Accepted loads in AHD cell--slot pair $a$ and locker $\ell$ \\
$C_\tau^{\mathrm{dist}}$, $C_\tau^{\mathrm{time}}$, $C_\tau^{\mathrm{fix}}$ & Distance, service-time, and fixed cost components & $C_k$ & Raw route-burden summary used in the value model \\
$\gamma$ & fulfillment-cost multiplier & $\phi(B_k)$, $z_j(B_k)$ & Aggregate basis and standardized basis component \\
$V_k(S_k)$ & Decision-state value function & $\widehat V^x(B_k)$ & Approximate post-decision value \\
$V_{k+1}^x(B)$ & Post-decision continuation value & $\theta_j$ & Linear value-function coefficient \\
$\Delta_k(S_k,o)$ & Exact opportunity cost of accepting option $o$ & $\pi^R$ & Rollout policy for value-function training \\
$\widehat\Delta_k(S_k,o)$ & Approximate opportunity cost from $\widehat V^x$ & $Y_i$; $\mathcal D$ & Continuation-profit label; labelled training set \\
$\widetilde c_{\eta,k}(o)$ & Adjusted opportunity cost used in pricing & $p^F(o)$ & Fine-static anchor price \\
$\eta\in[0,1]$ & ADP correction scale & $G_k^{\mathrm{TR}}(o)$; $\tau_h,\tau_\ell$ & Trust-region price grid; home and locker trust radii \\
\addlinespace
\midrule
\multicolumn{4}{@{}l}{\textit{Key abbreviations}}\\
\midrule
DOPMDO & Dynamic Offering and Pricing of Mixed Delivery Options & SVAP--ADP & State-Value Anchored Pricing via Approximate Dynamic Programming \\
AHD & Attended home delivery & OOH & Out-of-home delivery \\
LSP & Logistics service provider & MDP & Markov decision process \\
MNL & Multinomial logit & ADP & Approximate dynamic programming \\
VFA & Value-function approximation & TR & Trust region \\
\bottomrule
\end{tabularx}
\end{table}

\clearpage

\section{Aggregate Basis for the SVAP--ADP Value Approximation}
\label{app:vfa_basis_terms}

This appendix reports the aggregate basis used in the final SVAP--ADP value approximation. The basis is generated from the six operational dimensions introduced in Equation~\eqref{eq:core-aggregate-state}. The intercept is always included and is excluded from the ridge penalty. Let $D$ denote the booking-horizon length and $k$ the current epoch. Let $n_k^H$ and $n_k^L$ be the numbers of accepted AHD and OOH orders, respectively, and let $n_k=n_k^H+n_k^L$ be the total number of accepted orders. Let $A_k^L$ be the set of active lockers, $m_{\ell k}$ the number of parcels assigned to locker $\ell$, and $B_k^L=\{\ell\in A_k^L:m_{\ell k}=1\}$ the set of singleton lockers. Let $H_k$ denote the number of active home-delivery cells, and let $\bar q_k^H$ and $q_k^{H,\max}$ denote the implemented mean and maximum AHD route-count summaries over the active home-delivery aggregation. Let $C_k$ denote the accumulated raw distance-cost burden.

Table~\ref{tab:vfa_basis_terms} lists the 17 non-intercept basis components used in the deployed value model. The booking-progress block should be interpreted as a low-order time basis for the continuation-value level rather than as five separate state concepts. Under a fixed horizon, several of these terms are affine transformations of the same epoch variable, but they are retained because they are part of the locked final feature subset. Option dependence enters through successor-state evaluation: for each candidate option, SVAP--ADP recomputes the same aggregate basis after rejection and after hypothetical acceptance, and the difference in predicted state value defines the opportunity-cost adjustment.

Note that the components are deterministic basis transformations, not mutually independent explanatory variables. Some terms are correlated or derived from common operational quantities, but they are retained because the value model is used for prediction rather than coefficient interpretation. Validation experiments with pruned bases led to weaker value-approximation accuracy and policy performance, so the 17-component basis is the smallest retained set that consistently delivered stable out-of-sample results.

\begin{table}[!htbp]
\centering
\caption{Final aggregate basis used in the SVAP--ADP value approximation.}
\label{tab:vfa_basis_terms}
\begin{tabularx}{\textwidth}{p{0.23\textwidth}Y}
\toprule
Operational dimension & Retained basis components \\
\midrule
Booking progress & Current epoch $k$; fraction of horizon elapsed $k/D$; remaining number of epochs $D-k$; remaining-horizon fraction $(D-k)/D$; squared remaining-horizon fraction $((D-k)/D)^2$. \\
Committed load & Accepted AHD orders $n_k^H$; accepted OOH orders $n_k^L$; total accepted orders $n_k$. \\
Modal balance & AHD share $n_k^H/\max\{1,n_k\}$; OOH share $n_k^L/\max\{1,n_k\}$. \\
OOH consolidation & Active lockers $|A_k^L|$; orders per active locker $n_k^L/\max\{1,|A_k^L|\}$; singleton-locker share $|B_k^L|/\max\{1,|A_k^L|\}$. \\
AHD dispersion & Active home-delivery cells $H_k$; mean accepted AHD load over active home routes or time-window routes $\bar{q}_k^H$; maximum accepted AHD load over active home routes or time-window routes $q_k^{H,\max}$. \\
fulfillment burden & Raw distance-cost burden per accepted order $C_k/\max\{1,n_k\}$. \\
\bottomrule
\end{tabularx}
\end{table}

\section{Offline Training Algorithm for SVAP--ADP}
\label{app:offline_training}

Algorithm~\ref{alg:svap_offline_training} summarizes the offline rollout procedure used to train the SVAP--ADP post-decision value approximation. For each training seed, we simulate a reference trajectory under the rollout policy $\pi^R$, sample post-decision booking states, clone each sampled state, and complete the remaining horizon under the same policy. The realized continuation profit from the cloned rollout becomes the supervised label, computed as in Eq.~\eqref{eq:rollout-label}.

\begin{algorithm}[t]
\caption{Offline rollout training of the SVAP--ADP value approximation}
\label{alg:svap_offline_training}
\begin{algorithmic}[1]
\Require Rollout policy $\pi^R$; training seeds $\mathcal{S}^{\mathrm{train}}$; state-sampling rule; label budget $N_{\max}$; aggregate basis $\phi(\cdot)$; ridge parameter $\lambda$.
\Ensure Trained post-decision value model $\widehat V^x$.
\State Initialize the training set $\mathcal{D}\gets\emptyset$.
\For{each seed $s\in\mathcal{S}^{\mathrm{train}}$ while $|\mathcal{D}|<N_{\max}$}
    \State Generate the booking instance and random streams for seed $s$.
    \State Simulate one trajectory under $\pi^R$ and sample post-decision states $\mathcal{I}_s$.
    \For{each sampled state $B\in\mathcal{I}_s$ while $|\mathcal{D}|<N_{\max}$}
        \State Clone $B$ and roll out to the terminal state under $\pi^R$.
        \State Compute the continuation label $y(B)$ using Eq.~\eqref{eq:rollout-label}.
        \State Add $\bigl(\phi(B),y(B)\bigr)$ to $\mathcal{D}$.
    \EndFor
\EndFor
\State Standardize the non-intercept basis components.
\State Fit the ridge regression value model and store the learned coefficients.
\end{algorithmic}
\end{algorithm}

The label in Eq.~\eqref{eq:rollout-label} is future-only: it removes revenue already earned before the sampled state and retains the terminal fulfillment cost still to be incurred from that state. Thus, the regression target is aligned with the post-decision continuation value $V^x(B)$ rather than with total episode profit. The budget $N_{\max}$ limits the number of cloned rollouts, which are the main offline computational cost. After standardization, ridge regression stabilizes estimation under correlated aggregate basis components. The fitted model $\widehat V^x$ is then fixed and used online to compute option-level value differences for SVAP--ADP pricing.

\section{Seattle instance parameters}
\label{app:parameters}

Table~\ref{tab:seattle_params} summarizes the key parameters used in the Seattle case study. The first group gives the calibrated MNL choice model. The home-delivery utilities before price are set to $(5.50,5.45,5.35)$ for the three home options, while the outside-option utility is normalized to zero. The price coefficient $\beta^p=-0.30$ captures customer sensitivity to delivery charges, so higher prices reduce the probability of selection. Locker alternatives have a baseline utility $\mu^\ell=0.80$ and a distance-sensitivity parameter $\beta^d=0.55$ per kilometer, which controls how strongly customers are discouraged by walking or traveling farther to a pickup point. On the revenue side, each accepted order generates a base revenue of $r=20$, with home prices chosen from $[0,16]$ and locker prices from $[-2,12]$, allowing the policy to use both home-delivery surcharges and limited locker discounts.

The remaining parameter groups define the operational environment and the dynamic-pricing implementation. Fulfillment costs are scaled by $\gamma=3$, with distance cost rate $0.0045$, service-time components of $1.5$ for home stops and $0.3$ for locker stops, and a fixed route or vehicle cost of $8.0$; final routes are settled by the Clarke--Wright savings algorithm at the end of each episode. The candidate-menu construction starts from a locker pool of size $12$, always retains the nearest accessible locker, and uses active and forecast-inactive layer caps of $(6,6)$ before offering up to $6$ locker options to the customer. Finally, the SVAP--ADP configuration uses an opportunity-cost scale $\eta=1.0$, clipped adjusted costs in $[0.01,30.0]$, and trust-region radii $(\tau_h,\tau_\ell)=(1.0,2.0)$ for home and locker price adjustments.

\begin{table}[!htbp]
\centering
\caption{Calibrated parameters for the Seattle case study.}
\label{tab:seattle_params}
\small
\begin{tabular}{lll}
\toprule
\multicolumn{3}{l}{\textit{MNL choice model}} \\
\midrule
$\mu_M,\;\mu_A,\;\mu_E$ & $5.50,\;5.45,\;5.35$ & Home utilities before price \\
$\mu^\ell$ & $0.80$ & Locker baseline utility \\
$\beta^p$ & $-0.30$ & Price sensitivity \\
$\beta^d$ & $0.55$ & Locker-distance sensitivity per km \\
$u_0$ & $0$ & Outside utility \\
\addlinespace
\multicolumn{3}{l}{\textit{Fulfillment-cost components}} \\
\midrule
$\gamma$ & $3.0$ & Fulfillment-cost post-multiplier \\
Distance rate & $0.0045$ & Per unit distance, $C^{\text{dist}}$ \\
Service rate & $(1.5,0.3)$ & Per home or locker stop, $C^{\text{Time}}$ \\
Fixed route cost & $8.0$ & Per dispatched route and vehicle, $C^{\text{fix}}$ \\
Routing backend & Clarke--Wright & End-of-episode route settlement \\
\addlinespace
\multicolumn{3}{l}{\textit{Revenue and pricing}} \\
\midrule
$r$ & $20$ & Base revenue per accepted order \\
$\pi_{out}$ & $5$ & Opt-out penalty per no-purchase choice \\
$\mathcal{G}_h$ & $[0, 16]$ & Home price grid \\
$\mathcal{G}_\ell$ & $[-2, 12]$ & Locker price grid \\
\addlinespace
\multicolumn{3}{l}{\textit{Candidate menu}} \\
\midrule
\(K^H\) & 2 & Feasible time slots retained by common insertion-cost screen\\
$K^{\text{pool}}$ & $12$ & Total locker candidate-pool size \\
$K^{\text{near}}$ & $1$ & Nearest accessible locker retained \\
$(K^{\text{act}},\;K^{\text{inact}})$ & $(6,\;6)$ & Active and forecast-inactive layer caps \\
$K^{\text{offer}}$ & $6$ & Offered locker cap \\
\addlinespace
\multicolumn{3}{l}{\textit{SVAP--ADP parameters}} \\
\midrule
Value basis & $17+1$ & Seventeen aggregate basis components plus intercept \\
$\eta$ & $1.0$ & ADP opportunity-cost scale \\
Cost clipping & $[0.01,\;30.0]$ & Adjusted opportunity-cost bounds \\
$(\tau_h,\;\tau_\ell)$ & $(1.0,\;2.0)$ & Home and locker trust-region radii \\
\bottomrule
\end{tabular}
\end{table}

\end{document}